\pdfoutput=1
\RequirePackage{ifpdf}
\ifpdf 
\documentclass[pdftex]{sigma}
\else
\documentclass{sigma}
\fi

\numberwithin{equation}{section}

\newtheorem{Theorem}{Theorem}[section]

\newtheorem{Lemma}[Theorem]{Lemma}
\newtheorem{Proposition}[Theorem]{Proposition}

\newtheorem{theoremA}{Theorem}

 { \theoremstyle{definition}
\newtheorem{Definition}[Theorem]{Definition}

\newtheorem{Example}[Theorem]{Example}
\newtheorem{Remark}[Theorem]{Remark} }

\usepackage{multirow}
\usepackage{amscd}
\usepackage{mathrsfs}
\usepackage{pgf,tikz}

\DeclareMathOperator{\R}{\mathbb R}

\begin{document}

\newcommand{\arXivNumber}{2309.14697}

\renewcommand{\PaperNumber}{011}

\FirstPageHeading

\ShortArticleName{On Invariants of Constant $p$-Mean Curvature Surfaces in the Heisenberg Group $H_{1}$}

\ArticleName{On Invariants of Constant $\boldsymbol{p}$-Mean Curvature Surfaces\\ in the Heisenberg Group $\boldsymbol{H_{1}}$}

\Author{Hung-Lin CHIU~$^{\rm ab}$, Sin-Hua LAI~$^{\rm c}$ and Hsiao-Fan LIU~$^{\rm de}$}

\AuthorNameForHeading{H.-L.~Chiu, S.-H.~Lai and H.-F.~Liu}

\Address{$^{\rm a)}$~Department of Mathematics, National Tsing Hua University, Hsinchu, Taiwan}
\EmailD{\href{mailto:hlchiu@math.nthu.edu.tw}{hlchiu@math.nthu.edu.tw}}
\Address{$^{\rm b)}$~National Center for Theoretical Sciences, Taipei, Taiwan}

\Address{$^{\rm c)}$~Fundamental Education Center, National Chin-Yi University of Technology,\\
\hphantom{$^{\rm c)}$}~Taichung, Taiwan}
\EmailD{\href{mailto:shlai@ncut.edu.tw}{shlai@ncut.edu.tw}}

\Address{$^{\rm d)}$~Department of Applied Mathematics and Data Science, Tamkang University,\\
\hphantom{$^{\rm d)}$}~New Taipei City, Taiwan}
\EmailD{\href{mailto:hfliu@mail.tku.edu.tw}{hfliu@mail.tku.edu.tw}}

\Address{$^{\rm e)}$~Department of Applied Mathematics, National Chung Hsing University, Taiwan}

\ArticleDates{Received April 15, 2024, in final form February 04, 2025; Published online February 18, 2025}

\Abstract{One primary objective in submanifold geometry is to discover fascinating and significant classical examples of $H_1$. In this paper which relies on the theory we established in [\textit{Adv. Math.} \textbf{405} (2022), 08514, 50~pages, arXiv:2101.11780] and utilizing the approach we provided for constructing constant $p$-mean curvature surfaces, we have identified intriguing examples of such surfaces. Notably, we present a complete description of rotationally invariant surfaces of constant $p$-mean curvature and shed light on the geometric interpretation of the energy $E$ with a lower bound.}


\Classification{53A10; 53C42; 53C22; 34A26}

\section{Introduction}
This article is an extension of the previous paper \cite{ChiuH/LiuH:2022}, in which we studied the constant $p$-mean curvature surfaces in the Heisenberg group $H_{1}$. In \cite{ChiuH/LiuH:2022}, we focused on the foundation of the theory and paid more attention to the investigation of $p$-minimal surfaces. However, in the present article, instead of theory, we mainly focus on the examples, including an approach to construct constant $p$-mean curvature surfaces.

Recall that the Heisenberg group $H_1$ is the space $\mathbb{R}^3$ with the associated group multiplication%
\[(x_1,y_1,z_1)\circ (x_2,y_2,z_2)
=(x_1+x_2,y_1+y_2, z_1+z_2+y_1x_2-x_1y_2),\]
which is a $3$-dimensional Lie group. The space of all
left-invariant vector fields is spanned by the following three vector fields:
\[\mathring{e}_1=\frac{\partial}{\partial x}+y\frac{\partial}{\partial z},\qquad
\mathring{e}_2=\frac{\partial}{\partial y}-x\frac{\partial}{\partial z}
\qquad\text{and}\qquad T=\frac{\partial}{\partial z}.\]
The Heisenberg dilation (scaling) by the factor $\delta>0$ is the map $D_\delta\colon H_1\rightarrow H_1$ defined by~${D_\delta(x,y,z)=\bigl(\delta x,\delta y,\delta^2 z\bigr)}$ for any $(x,y,z)\in H_1$ (see \cite{Dragomir/Tomassini:2006}).

The standard contact bundle on $H_1$ is the subbundle $\xi$ of the tangent bundle
$TH_1$ spanned by $\mathring{e}_1$ and $\mathring{e}_2$. It is also defined to be the kernel of
the contact form
$\Theta={\rm d}z+x{\rm d}y-y{\rm d}x$.
The CR structure on $H_1$ is the endomorphism $J\colon \xi\to \xi$ defined by
$J(\mathring{e}_1)=\mathring{e}_2$ and $J(\mathring{e}_2)=-\mathring{e}_1$.
One can view $H_1$ as a pseudo-hermitian manifold
with $(J,\Theta)$ as the standard pseudo-hermitian structure. There is a naturally associated connection $\nabla$ if we regard all these left-invariant vector fields~$\mathring{e}_1$,~$\mathring{e}_2$, and $T$ as parallel vector fields.
A naturally associated metric on $H_{1}$ is the adapted metric $g_{\Theta}$, which is defined by $g_{\Theta}={\rm d}\Theta(\cdot,J\cdot)+\Theta^{2}$. It is equivalent to defining the metric regarding $\mathring{e}_1$, $\mathring{e}_2$, and $T$ as an orthonormal frame field. We sometimes use $\langle \cdot,\cdot\rangle$ to denote the adapted metric. In this paper, we use the adapted metric to measure the lengths, angles of vectors, and so on.

Suppose $\Sigma$ is a surface in the Heisenberg group $H_{1}$. There is a one-form $I$ on $\Sigma$ induced from the adapted metric $g_{\Theta}$. This induced metric is defined on the whole surface $\Sigma$ and is called the first fundamental form of $\Sigma$. The intersection $T\Sigma\cap\xi$ is integrated to be a singular foliation on~$\Sigma$ called the characteristic foliation. Each leaf is called a characteristic curve. A point $p\in\Sigma$ is called a singular point if the tangent plane $T_{p}\Sigma$ coincides with the contact plane $\xi_{p}$; otherwise, $p$ is called a regular (or non-singular) point. Generically, a point $p\in\Sigma$ is a regular point, and the set of all regular points is called the regular part of $\Sigma$. In this paper, we always assume that the surface $\Sigma$ is of class $C^{2}$, but of class $C^{\infty}$ on the regular part.
On the regular part, we can choose a unit vector field $e_{1}$ such that $e_{1}$ defines the characteristic foliation. The vector $e_{1}$ is determined up to a sign. Let $e_{2}=Je_{1}$. Then $\{e_{1},e_{2}\}$ forms an orthonormal frame field of the contact bundle $\xi$. We usually call the vector field $e_{2}$ a horizontal normal vector field. Then the $p$-mean curvature $H$ of the surface $\Sigma$ is defined by
$
\nabla_{e_{1}}e_{2}=-He_{1}$.
The $p$-mean curvature $H$ is only defined on the regular part of $\Sigma$. If $H=c$, which is a~constant on the whole regular part, we call the surface a constant $p$-mean curvature surface. In particular, if $c=0$, it is a~$p$-minimal surface. There also exists a function $\alpha$ defined on the regular part such that $\alpha e_{2}+T$ is tangent to the surface~$\Sigma$. We call this function the $\alpha$-function of~$\Sigma$. It is uniquely determined up to a sign, which depends on the choice of the characteristic direction~$e_{1}$. Define~${\hat{e}_{1}=e_{1}}$ and \smash{$\hat{e}_{2}=\frac{\alpha e_{2}+T}{\sqrt{1+\alpha^2}}$}, then $\{\hat{e}_{1},\hat{e}_{2}\}$ forms an orthonormal frame field of the tangent bundle $T\Sigma$. Notice that $\hat{e}_{2}$ is uniquely determined and independent of the choice of the characteristic direction $e_{1}$. In \cite{ChiuH/HuangY/LaiS:2017,ChiuH/LaiS:2015}, it was shown that these three invariants,
$I$, $e_{1}$, and $\alpha$,
form a complete set of invariants for constant $p$-mean curvature surfaces with $H=c$ in $H_{1}$. Namely, for any two surfaces with the same constant $p$-mean curvature having the same $I$, $\alpha$, $e_1$, they are differed only by a Heisenberg symmetry. In particular, if $\Sigma\subset H_{1}$ is a \emph{constant} $p$-mean curvature surface with $H=c$, then in terms of a compatible coordinate system $(U; x,y)$, which means
$e_{1}=\frac{\partial}{\partial x}$,
 the integrability condition (see \cite{ChiuH/LiuH:2022}) is reduced to
\begin{gather}
-a_{x}+a\frac{b_{x}}{b}=\frac{c\alpha}{\bigl(1+\alpha^2\bigr)^{1/2}},\qquad
-\frac{b_{x}}{b}=2\alpha+\frac{\alpha\alpha_{x}}{1+\alpha^2},\nonumber\\
\alpha_{xx}+6\alpha\alpha_{x}+4\alpha^{3}+c^{2}\alpha=0,\label{intcon}
\end{gather}
where the two functions $a$ and $b$ are a representation of the first fundamental form $I$ in the following sense that they describe the vector field
\begin{equation*}
\hat{e}_{2}=a(x, y)\frac{\partial}{\partial x}+b(x, y)\frac{\partial}{\partial y}.
\end{equation*}
In other words, there exists the $\alpha$ satisfying the \emph{Codazzi-like} equation
\begin{equation}\label{Codeq}
\alpha_{xx}+6\alpha\alpha_{x}+4\alpha^{3}+c^{2}\alpha=0,
\end{equation}
which is a nonlinear ordinary differential equation. In \cite{ChiuH/LiuH:2022}, we normalized $a$ and $b$ such that they can be uniquely determined by the function $\alpha$, and hence we obtained the result that the existence of a constant $p$-mean curvature surface (without singular points) is equivalent to the existence of a solution to a nonlinear second-order ODE \eqref{Codeq}, which is a kind of \emph{Li\'{e}nard equations} (cf.\ \cite{Polyanin/Zaitsev:2003}). They are one-to-one correspondences in some sense. For a detailed description, see~\cite[Theorems~1.1, 1.3 and~6.3]{ChiuH/LiuH:2022}. This result tells us that the investigation of the geometry of constant $p$-mean curvature surfaces in $H_{1}$ is equal to the study of the solution of the equation~\eqref{Codeq}. More specifically, we obtained a complete set of solutions (see \cite[Theorems~1.2 and~1.4]{ChiuH/LiuH:2022} or Theorem~\ref{consta1}) and used the types of the solutions to the equation to characterize the constant $p$-mean curvature surfaces as several classes, which are \emph{vertical}, \emph{special type I}, \emph{special type II} and \emph{general type} (see \cite[Definitions 5.1 and 5.2]{ChiuH/LiuH:2022} for $p$-minimal cases and see Definitions~\ref{cladef} and~\ref{cladef1} for the cases with $c\neq 0$ in the present article). After the process of normalization, we obtained a complete set of invariants from the normal form of the $\alpha$-function. It is worth of our mention that these invariants in some sense measure how different a constant $p$-mean curvature surface is from the model case, which is the horizontal plane in the $p$-minimal case, and the Pansu sphere in the case $c>0$.

We first study rotationally invariant surfaces in $H_1$ with constant $p$-mean curvature $H=c$ using the Codazzi-like equation \eqref{Codeq}. In \cite{Ritore/Rosales:2006}, M.~Ritor\'{e} and C.~Rosales made an investigation on such kinds of surfaces by a first-order ODE system. In the present paper, we shall study them again from the point of view of our theory established in the previous paper \cite{ChiuH/LiuH:2022} and the present one.
Let $\Sigma(s,\theta)$ be a rotationally invariant surface in $H_1$ with $H=c$, generated by the curve~${\gamma(s)=(x(s),0,t(s))}$, $x(s)\geq 0$, on the $xt$-plane,
that is, $\Sigma$ is parametrized by
\begin{equation*}
\Sigma(s,\theta)=(x(s)\cos \theta,x(s)\sin \theta,t(s)),
\end{equation*}
where $x'^{2}+t'^{2}=1$. Here $'$ means taking a derivative with respect to $s$. Recall the energy
\begin{equation}\label{energy1}
E=\frac{xt^{\prime }}{\sqrt{x^2x^{\prime 2}+t^{\prime 2}}}+\lambda x^2,
\end{equation}
which was introduced in \cite{Ritore/Rosales:2006} and was shown to be a constant. Here $2\lambda=c$ and notice that our $p$-mean curvature differs from the one defined in \cite{Ritore/Rosales:2006} by a sign. Hence, we have Theorems~\ref{TheoremA} and~\ref{TheoremB} as follows.

\begin{theoremA}\label{TheoremA}
A curve $\gamma=(x,t)$ is the generating curve of a rotationally invariant surface $\Sigma$ in~$H_{1}$ with $H=c\neq 0$
if and only if $\gamma=(x,t)$ is defined by
$x^{2}=\frac{k}{c^2}+r\cos{(c\tilde{s})}$, $t=-\frac{\tilde{s}}{c}-\frac{r}{2}\sin{(c\tilde{s})}$, up to a constant,
for some horizontal arc-length parameter $\tilde{s}$ and some $k, r\in\R$ such that
\begin{equation*}
k\geq 1\qquad \textrm{and}\qquad r=\frac{2}{c^2}\sqrt{k-1}.
\end{equation*} In addition, we have
$k=2cE+2$. If $r=0$, then $\Sigma$ is a cylinder. If $r\neq 0$, then, in terms of normal coordinates $\bigl(\bar s,\bar\theta\bigr)$, the two invariants for $\Sigma$ are
\begin{alignat*}{3}
&\zeta_{1}(\bar{\theta})=-\frac{2E\bar\theta}{cr},\qquad && \textrm{up to a constant, which is linear on $\bar\theta$}, &\\
& \zeta_{2}(\bar{\theta})=-\frac{2cE+2}{c^{2}r},\qquad&& \textrm{which is a constant}.&
\end{alignat*}
\end{theoremA}
\begin{theoremA}\label{TheoremB}
A curve $\gamma=(x,t)$ is the generating curve of a rotationally invariant $p$-minimal surface $\Sigma$ in $H_{1}$ if and only if either $t$ is a constant, and hence $\Sigma$ is a part of the horizontal plane, or $\gamma=(x,t)$ is defined by
$
x^{2}=\tilde{s}^{2}+c_2$, $
t=m\tilde{s}$, up to a constant,
for some horizontal arc-length parameter $\tilde{s}$ and some $c_2, m\in\R,\ m\neq 0$. In addition, we have~${
E=m}$.
In terms of normal coordinates $(\bar s,\bar\theta)$, the two invariants for $\Sigma$ are
\begin{alignat*}{3}
& \zeta_{1}(\bar{\theta})=E\bar\theta,\qquad && \textrm{up to a constant, which is linear on $\bar\theta$},& \\
&\zeta_{2}(\bar{\theta})=c_{2},\qquad && \textrm{which is a constant}.&
\end{alignat*}
\end{theoremA}
For more interesting examples, in Section \ref{ccpms}, we provide an approach to construct a constant $p$-mean curvature surface. This approach is an analog of the one we performed in the previous paper \cite{ChiuH/LiuH:2022} for $p$-minimal surfaces. Actually, in \cite{ChiuH/LiuH:2022}, we deformed the horizontal plane along a~curve~${\mathcal{C}(\theta)=(x_{1}(\theta), x_{2}(\theta), x_{3}(\theta))}$ to obtain a $p$-minimal surface.
More specifically, in \cite[Section 9]{ChiuH/LiuH:2022}, depending on a parametrized curve $\mathcal{C}(\theta)=(x_{1}(\theta), x_{2}(\theta), x_{3}(\theta))$ for $\theta\in\R$, we deformed the graph $u=0$ to obtain a $p$-minimal surface parametrized by
 \begin{equation*}
 Y(r,\theta)=(x_{1}(\theta)+r\cos{\theta},x_{2}(\theta)+r\sin{\theta},x_{3}(\theta)+rx_{2}(\theta)\cos{\theta}-rx_{1}(\theta)\sin{\theta}),
 \end{equation*}
for $r\in\R$. It is easy to check that $Y$ is \emph{an immersion} if and only if either $\Theta\bigl(\mathcal{C}'(\theta)\bigr)-\bigl(x_{2}'(\theta)\cos{\theta}\smash{-x_{1}'(\theta)\sin{\theta}\bigr)^{2}\neq 0}$ or $r+\bigl(x_{2}'(\theta)\cos{\theta}-x_{1}'(\theta)\sin{\theta}\bigr)\neq 0$ for all $\theta$.
In particular, the surface $Y$ defines a $p$-minimal surface of \emph{special type I} if the curve $\mathcal{C}$ satisfies
\begin{equation*}
x_{3}'(\theta)+x_{1}(\theta)x_{2}'(\theta)-x_{2}(\theta)x_{1}'(\theta)-\bigl(x_{2}'(\theta)\cos{\theta}-x_{1}'(\theta)\sin{\theta}\bigr)^2=0,
\end{equation*}
for all $\theta$.
In addition, the corresponding $\zeta_{1}$-invariant \cite[formula~(9.9)]{ChiuH/LiuH:2022} reads
\begin{equation}\label{713i}
\zeta_{1}(\theta)=x_{2}'(\theta)\cos{\theta}-x_{1}'(\theta)\sin{\theta}-\int \bigl[x_{1}'(\theta)\cos{\theta}+x_{2}'(\theta)\sin{\theta}\bigr]{\rm d}\theta.
\end{equation}

Similarly, the surface $Y$ defines a $p$-minimal surface of \emph{general type} if the curve $\mathcal{C}$ satisfies
\begin{equation*}
x_{3}'(\theta)+x_{1}(\theta)x_{2}'(\theta)-x_{2}(\theta)x_{1}'(\theta)-\bigl(x_{2}'(\theta)\cos{\theta}-x_{1}'(\theta)\sin{\theta}\bigr)^2\neq 0,
\end{equation*}
for all $\theta$. In addition, the corresponding $\zeta_{1}$- and $\zeta_{2}$-invariant read
\begin{gather}
\zeta_{1}(\theta)=x_{2}'(\theta)\cos{\theta}-x_{1}'(\theta)\sin{\theta}-\int \bigl[x_{1}'(\theta)\cos{\theta}+x_{2}'(\theta)\sin{\theta}\bigr] {\rm d}\theta,\nonumber\\
\zeta_{2}(\theta)=x_{3}'(\theta)+x_{1}(\theta)x_{2}'(\theta)-x_{2}(\theta)x_{1}'(\theta)-\bigl(x_{2}'(\theta)\cos{\theta}-x_{1}'(\theta)\sin{\theta}\bigr)^2\label{714i}
.
\end{gather}

In Section \ref{ccpms}, we construct a constant $p$-mean curvature surface by perturbing the Pansu sphere along a given curve $\mathcal{C}(\theta)$. In Section \ref{pansu1}, we see that the Pansu sphere \eqref{pansusphere} can be parametrized by
\begin{gather*} 
X(s,\theta)= \left(
\begin{matrix}
\cos \theta & -\sin \theta & 0 \\
\sin \theta & \cos \theta & 0 \\
0 & 0 & 1
\end{matrix}
\right)\left(
\begin{matrix}
x(s) \\
y(s) \\
t(s)
\end{matrix}
\right)=\left(
\begin{matrix}
x(s)\cos \theta-y(s)\sin \theta \\
x(s)\sin \theta+y(s)\cos \theta \\
t(s)
\end{matrix}
\right),
\end{gather*}
where
\begin{gather*} 
x(s) = \frac{1}{2\lambda}\sin(2\lambda s),\qquad
y(s) = -\frac{1}{2\lambda}\cos(2\lambda s)+ \frac{1}{2\lambda},\\
t(s) = \frac{1}{4\lambda^2}\sin(2\lambda s)- \frac{1}{2\lambda}s+\frac{\pi}{4\lambda^2},
\end{gather*}
with $X(0,\theta)=\bigl(0,0,\frac{\pi}{4\lambda^2}\bigr)$ and $X\bigl(\frac{\pi}{\lambda},\theta\bigr)=\bigl(0,0,-\frac{\pi}{4\lambda^2}\bigr)$ as the North pole and South pole, respectively.
We deform it along $\mathcal{C}(\theta)=(x_{1}(\theta), x_{2}(\theta), x_{3}(\theta))$ to obtain a constant $p$-mean curvature surface $Y(s,\theta)$ as follows
\begin{align*} 
Y(s,\theta)={}&(x_{1}(\theta)+(x(s)\cos{\theta}-y(s)\sin{\theta}),x_{2}(\theta)+(x(s)\sin{\theta}+y(s)\cos{\theta}), \\
&x_{3}(\theta)+t(s)+x_{2}(\theta)(x(s)\cos{\theta}-y(s)\sin{\theta})
-x_{1}(\theta)(x(s)\sin{\theta}+y(s)\cos{\theta})).
\end{align*}
We also give a condition for $Y$ to be an immersion. The coordinate system $(s,\theta)$ for $Y$
is a~compatible one. We have (see \eqref{forinv})
\begin{equation}\label{forinvi}
\alpha =\lambda \frac{A(\theta )\cos {2\lambda s}+\left( \frac{1}{2\lambda }-B(\theta )\right) \sin {2\lambda s}}{\left( B(\theta )-\frac{1}{2\lambda }\right) \cos {2\lambda s}+A(\theta )\sin {2\lambda s}+D(\theta )},
\end{equation}
where
\begin{gather*} 
A(\theta)=x^{\prime }_{2}(\theta)\cos{\theta}-x^{\prime }_{1}(\theta)\sin{\theta},\qquad B(\theta)=x^{\prime }_{2}(\theta)\sin{\theta}+x^{\prime }_{1}(\theta)\cos{\theta}, \\
D(\theta)=\lambda \Theta(\mathcal{C}^{\prime }(\theta))+\left(\frac{1}{2\lambda}-B(\theta)\right),\qquad V(\theta)=\left( A(\theta ),\frac{1}{2\lambda }-B(\theta )\right).
\end{gather*}
It is obvious that $V(\theta )=0$ implies $\alpha=0$, and hence $Y$ is a cylinder. For nonzero $V(\theta )$, we then define
\[\left\Vert V(\theta)\right \Vert =\sqrt{\left[ A( \theta ) \right] ^{2}+
\left[ \frac{1}{2\lambda }-B(\theta )\right] ^{2}}\]
and write
\[
G(\theta )=\frac{D(\theta )}{ \Vert V(\theta) \Vert },\qquad \frac{V(\theta)}{ \Vert V(\theta) \Vert }=(\sin {\zeta (\theta )},\cos {\zeta
(\theta )}),
\]
for some function $\zeta {(\theta )}$.
From \eqref{forinvi}, we thus have
\[\alpha=\lambda \frac{\sin {\bigl(2\lambda s+\zeta {(\theta )}\bigr)}}{G(\theta
)-\cos {\bigl(2\lambda s+\zeta {(\theta )}\bigr)}}.\]
Finally, we normalize the $\alpha$-function to obtain the two invariants for $Y$, stated in Theorem \ref{construction}. We list it here as the third main result of the present paper.

\begin{theoremA}\label{TheoremC}
If $V=0$, then $\alpha=0$. If $V\neq 0$, we consider the new coordinates
$
\tilde{s}=s+\Gamma(\theta)$, $ \tilde{\theta}=\Psi( \theta)$,
where
\[
\Gamma(\theta)=\frac{\theta }{2\lambda}-\int D( \theta){\rm d}\theta,\qquad
\Psi(\theta)=2\lambda \int \left \Vert V(\theta)\right \Vert {\rm d}\theta ,
\]
then the coordinate system $\bigl(\tilde{s},\tilde{\theta}\bigr)$ is normal. In terms of the
normal coordinates, the invariants of $Y$ are given by
\begin{equation}\label{invariantsi}
\zeta _{1}\bigl( \tilde{\theta}\bigr) = \zeta (\theta) -2\lambda \Gamma (\theta), \qquad
\zeta _{2}\bigl( \tilde{\theta}\bigr) = G(\theta),
\end{equation}
where $\theta= \Psi ^{-1}\bigl( \tilde{\theta}\bigr)$.
\end{theoremA}

In Section \ref{examples}, we use formula \eqref{713i}, \eqref{714i} and \eqref{invariantsi} to construct various examples of surfaces of constant $p$-mean curvature including degenerate p-minimal surfaces of special type I.

\section{Solutions to the Codazzi-like equation}

The Codazzi-like equation for a surface in $H_1$ with constant $p$-mean
curvature $H=c>0$ is
\begin{equation} \label{ceq}
\alpha_{xx}+6\alpha \alpha_{x}+4\alpha^{3}+c^2\alpha=0.
\end{equation}

\begin{Theorem}[\cite{ChiuH/LiuH:2022}]
\label{consta1} Besides the following three special solutions to \eqref{ceq},
\[
\alpha(x)=0,\qquad -\frac{c}{2}\tan(cx+cK_1),\qquad \alpha(x)=-\frac{c}{2}\tan\left(\frac{c%
}{2}x-\frac{c}{2}K_2\right),
\]
we have the general solution to \eqref{ceq} of the form
\begin{equation*}
\alpha(x)=\frac{c}{2}\frac{\sin{\bigl(cx+c_{1}\bigr)}}{c_{2}-\cos{\bigl(cx+c_{1}%
\bigr)}},
\end{equation*}
which depends on constants $K_{1}$, $K_{2}$, $c_1$, and $c_2$.
\end{Theorem}

Note that all the solutions are a periodic function with $\alpha\bigl(x+\frac{2\pi%
}{c}\bigr)=\alpha(x)$ for all $x$. We give some remarks as follows.
\begin{enumerate}\itemsep=0pt
\item[(1)] In terms of the following identities
\begin{equation*}
-\tan{\left(\theta+\frac{\pi}{2}\right)}=\cot{\theta}=\frac{\sin{2\theta}}{1-\cos{2\theta}},\qquad
-\tan{2\theta}=-\frac{\sin{2\theta}}{\cos{2\theta}}=\frac{\sin{2\theta}}{0-\cos{2\theta}},
\end{equation*}
we see that the two nontrivial special solutions in Theorem \ref{consta1}
correspond to the general solution in Theorem \ref{consta1} with $c_{2}=0$ and $c_{2}=1$, respectively.
\item[(2)] From the following identity
\[
\frac{\sin{(\theta+\pi)}}{c_{2}-\cos{(\theta+\pi)}}=\frac{\sin{\theta}}{%
-c_{2}-\cos{\theta}},
\]
we can assume without loss of generality that $c_{2}\geq 0$ in the general
solution.
\end{enumerate}

Due to Theorem \ref{consta1}, we are able to use the types of the solutions
to \eqref{ceq} to classify the constant $p$-mean curvature surfaces into
several classes, which are vertical, special type I, special type II and
general type. In terms of compatible coordinates $(x,y)$, the function $\alpha(x,y)$ is a~solution to the \emph{Codazzi-like} equation \eqref{ceq}
for any given $y$. By Theorem \ref{consta1}, the function~$\alpha(x,y)$
hence has one of the following forms of special types
\[
0,\qquad \frac{c}{2}\frac{\sin(cx+c_1)}{0-\cos(cx+c_1)},\qquad
\frac{c}{2}\frac{\sin{(cx+c_{1})}}{1-\cos{\bigl(cx+c_{1})}},
\]
and general types
$
\frac{c}{2}\frac{\sin{(cx+c_{1})}}{c_{2}-\cos{(c(x+c_{1})}}$,
where, instead of constants, both $c_{1}$ and $c_{2}$ are now functions of $%
y $. Notice that it is convenient at some point to assume that $c_{2}(y)\geq0$ for all $y$. We now use the types of the
function $\alpha(x,y)$ to define the types of constant $p$-mean curvature
surfaces as follows.

\begin{Definition}
\label{cladef} Locally, we say that a constant $p$-mean curvature surface is
\begin{enumerate}\itemsep=0pt
\item[(1)] \emph{vertical} if $\alpha$ vanishes (i.e., $\alpha(x,y)=0$ for all $x,y$);
\item[(2)] of \emph{special type I} if $\alpha=\frac{c}{2}\frac{\sin{(cx+c_1(y))}}{1-\cos{(cx+c_1(y))}}$;
\item[(3)] of \emph{special type II} if $\alpha=\frac{c}{2}\frac{\sin(cx+c_1(y))}{0-\cos(cx+c_1(y))}$;
\item[(4)] of \emph{general type} if $\alpha=\frac{c}{2}\frac{\sin{(cx+c_1(y))}}{c_{2}(y)-\cos{(cx+c_1(y))}}$ with $c_{2}(y) \notin \{0,1\}$
for all $y$.
\end{enumerate}
\end{Definition}

We further divide constant $p$-mean curvature surfaces of general type into
three classes as follows.

\begin{Definition}
\label{cladef1} A constant $p$-mean curvature surface of \emph{general type%
} is
\begin{enumerate}\itemsep=0pt
\item[(1)] of \emph{type I} if $c_{2}(y)>1$ for all $y$;
\item[(2)] of \emph{type II} if $0<c_{2}(y)<1$ for all $y$, and $\frac{-c_{1}+\cos^{-1}c_2}{c}<x<\frac{2\pi-c_{1}-\cos^{-1}c_2}{c}$;
\item[(3)] of \emph{type III} if $0<c_{2}(y)<1$ for all $y$, and either $\frac{-c_{1}}{c}\leq x<\frac{-c_{1}+\cos^{-1}c_2}{c}$ or
$\frac{2\pi-c_{1}-\cos^{-1}c_2}{c}<x\leq \frac{-c_{1}+2\pi}{c}$,
\end{enumerate}
where $\cos^{-1}$ is the inverse of the function $\cos\colon [0,\pi]\rightarrow[
-1,1]$.
\end{Definition}

\begin{figure}[t]
\centering
\includegraphics[scale=0.5]{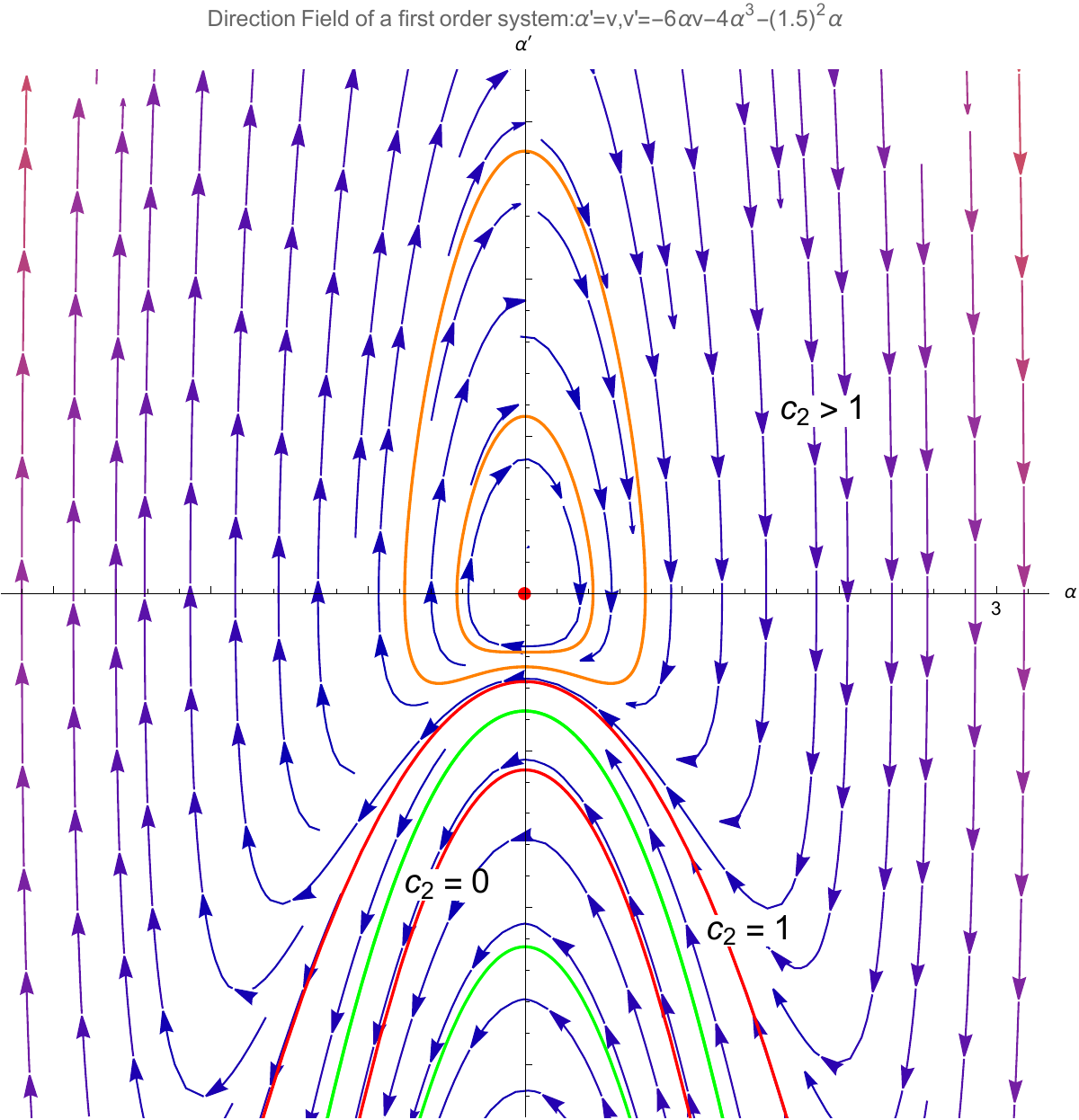}
\caption{Direction field for $c=1.5$.}\label{fig1}
\end{figure}

We notice that the \emph{type} is invariant under the action of a
Heisenberg rigid motion and the regular part of a constant $p$-mean
curvature surface $\Sigma \subset H_{1}$ is a union of these types of
surfaces. The corresponding paths of each type of $\alpha$ are shown on the
phase plane (see Figure \ref{fig1}). We express some basic facts as follows.

\begin{itemize}\itemsep=0pt
\item If $\alpha$ vanishes, then it is part of a vertical cylinder.

\item The two concave downward parabolas in red represent
\[
\alpha=\frac{c}{2}\frac{\sin{(cx+c_1)}}{1-\cos{(cx+c_1)}},\qquad \frac{c}{2}\frac{\sin{(cx+c_1)}}{0-\cos{(cx+c_1)}},
\]
respectively. The one for $\alpha=\frac{c}{2}\frac{\sin{(cx+c_1)}}{1-\cos{(cx+c_1)}}$ is above the one for $\alpha=\frac{c}{2}\frac{\sin{(cx+c_1)}}{0-\cos{(cx+c_1)}}$. For
surfaces of \emph{special type I}, we have that
\smash{$
\alpha\left(\frac{\pi-c_1}{c}\right)=0$}, \smash{$ \alpha^{\prime }\left(\frac{\pi-c_1}{c}\right)=-\frac{c^2}{4}
$}
and
\begin{equation*}
\alpha \rightarrow \begin{cases}
\infty, & \text{if}\ x\rightarrow\dfrac{-c_{1}}{c}\quad \text{from the right}, \vspace{1mm}\\
-\infty, & \text{if}\ x\rightarrow \dfrac{2\pi-c_1}{c} \quad \text{from the
left},
\end{cases}
\end{equation*}
and, for surfaces of \emph{special type II} in which $\alpha$ has period $\pi$, we have that
\[
\alpha\left(\frac{-c_{1}}{c}\right)=\alpha\left(\frac{\pi-c_1}{c}\right)=0,\qquad \alpha^{\prime
}\left(\frac{-c_{1}}{c}\right)=\alpha^{\prime }\left(\frac{\pi-c_1}{c}\right)=-\frac{c^2}{2}
\]
and
\begin{equation*}
\alpha \rightarrow \begin{cases}
\infty, & \text{if}\ x\rightarrow \dfrac{\pi-2c_1}{2c}\quad \text{from the
right}, \vspace{1mm}\\
-\infty, & \text{if}\ x\rightarrow \dfrac{3\pi-2c_1}{2c}\quad \text{from the
left}.
\end{cases}
\end{equation*}

\item The closed curves in orange on the phase plane correspond to the family of solutions
\[
\alpha(x)=\frac{c}{2}\frac{\sin{(cx+c_{1})}}{c_{2}-\cos{(cx+c_{1})}},
\]
where $c_1$, $c_2$ are constants and $c_2>1$, which are of \emph{type I}.
There exist zeros for $\alpha$-function at $x=\frac{-c_{1}}{c},\ \frac{\pi-c_1}{c}$,
at which we have that
\[
\alpha^{\prime }\left(\frac{-c_1}{c}\right)=\frac{1}{2}\frac{c^2}{c_{2}-1}>0,\qquad \alpha^{\prime
}\left(\frac{\pi-c_1}{c}\right)=-\frac{1}{2}\frac{c^2}{c_{2}+1}<0.
\]
There are no singular points for surfaces of \emph{type I}.

\item The curves in between the two red concave downward parabolas are
of \emph{type II}. The $\alpha$-function of \emph{type II} has a zero at
$x=\frac{\pi-c_1}{c}$, and
\[
\alpha^{\prime }\left(\frac{\pi-c_1}{c}\right)=-\frac{1}{2}\frac{c^2}{c_{2}+1}<0.
\]

For surfaces of \emph{type II}, it can be checked that
\begin{equation*}
\alpha \rightarrow \begin{cases}
\infty, & \text{if}\ x\rightarrow \dfrac{-c_{1}+\cos^{-1}c_2}{c}\quad \text{from the right}, \vspace{1mm}\\
-\infty, & \text{if}\ x\rightarrow \dfrac{-c_{1}+2\pi-\cos^{-1}c_2}{c}\quad
\text{from the left}.
\end{cases}
\end{equation*}

\item The curves beneath the lower concave downward parabola are of \emph{type III}. There exists a~zero for $\alpha$-function at $x=\frac{-c_{1}}{c}$, and
\[
\alpha^{\prime }\left(\frac{-c_1}{c}\right)=\frac{1}{2}\frac{c^2}{c_{2}-1}<0.
\]
For surfaces of \emph{type III}, we have
\begin{equation*}
\alpha \rightarrow\begin{cases}
-\infty, & \text{if}\ x\rightarrow \dfrac{-c_{1}+\cos^{-1}c_2}{c}\quad \text{from the left}, \vspace{1mm}\\
\infty, & \text{if}\ x\rightarrow \dfrac{-c_{1}+2\pi-\cos^{-1}c_2}{c}\quad
\text{from the right}.
\end{cases}
\end{equation*}
\end{itemize}

\subsection{The Pansu sphere}\label{pansu1}

\begin{Lemma}[Pansu sphere]
A Pansu sphere given in {\rm\cite{ChengJ/ChiuH/HwangJ/YangP:2018}} by
\begin{equation}\label{pansusphere}
f(z)=\frac{1}{2\lambda^2}\bigl(\lambda|z|\sqrt{1-\lambda^2|z|^2}+\cos^{-1}(\lambda|z|)\bigr),\qquad |z|\leq \frac{1}{\lambda},
\end{equation}
of constant $p$-mean curvature $c=2\lambda$ has its $\alpha$-function of special type I. In fact, we have
\begin{equation*}
\alpha=\frac{\lambda \sin(2\lambda s)}{(1-\cos(2\lambda s))},\qquad a=\frac{-\lambda}{\sqrt{1+\alpha^2}},\qquad b=\frac{2\lambda^{2}}{\sqrt{1+\alpha^2}(1-\cos{2\lambda s})}.
\end{equation*}
\end{Lemma}

\begin{proof}
We parametrize a Pansu sphere by
\begin{equation} \label{paps}
X(s,\theta)= \left(
\begin{matrix}
\cos \theta & -\sin \theta & 0 \\
\sin \theta & \cos \theta & 0 \\
0 & 0 & 1
\end{matrix}
\right)\left(
\begin{matrix}
x(s) \\
y(s) \\
t(s)
\end{matrix}
\right)=\left(
\begin{matrix}
x(s)\cos \theta-y(s)\sin \theta \\
x(s)\sin \theta+y(s)\cos \theta \\
t(s)
\end{matrix}
\right),
\end{equation}
where
\begin{gather}
x(s) = \frac{1}{2\lambda}\sin(2\lambda s),\qquad
y(s) = -\frac{1}{2\lambda}\cos(2\lambda s)+ \frac{1}{2\lambda},\nonumber\\
t(s) = \frac{1}{4\lambda^2}\sin(2\lambda s)- \frac{1}{2\lambda}s+\frac{
\pi}{4\lambda^2}, \label{ad}
\end{gather}
with $X(0,\theta)=\bigl(0,0,\frac{\pi}{4\lambda^2}\bigr)$ and $X\bigl(\frac{\pi}{\lambda},\theta\bigr)=\bigl(0,0,-\frac{\pi}{4\lambda^2}\bigr)$ as the North pole and South pole, respectively.
We then have $e_1:=X_s$, which means that $X(s,\theta)$
defines a compatible coordinate system. Moreover, $e_2=Je_1$, where
\begin{align*}
e_1&=X_s=\bigl(x^{\prime }(s)\cos \theta-y^{\prime }(s)\sin \theta,x^{\prime
}(s)\sin \theta+y^{\prime }(s)\cos \theta,t^{\prime }(s)\bigr) \\
&=\bigl(x^{\prime }(s)\cos \theta-y^{\prime }(s)\sin \theta\bigr) \mathring
e_1+\bigl(x^{\prime }(s)\sin \theta+y^{\prime }(s)\cos \theta\bigr)\mathring
e_2, \\
e_2&=-\bigl(x^{\prime }(s)\sin \theta+y^{\prime }(s)\cos \theta\bigr)\mathring
e_1+\bigl(x^{\prime }(s)\cos \theta-y^{\prime }(s)\sin \theta\bigr) \mathring e_2.
\end{align*}
We note that $\alpha$ is a function satisfying
\begin{equation} \label{aa}
\alpha e_2+T=\mathcal{A}X_s+\mathcal{B}X_\theta,
\end{equation}
for some functions $\mathcal{A}$ and $\mathcal{B}$. Direct calculation shows that
\begin{equation*}
X_\theta=(-x(s)\sin \theta-y(s)\cos \theta)\mathring e_1+(x(s)\cos \theta-y(s)\sin
\theta)\mathring e_2+(x^2(s)+y^2(s))T.
\end{equation*}
Therefore, \eqref{aa} implies
\begin{gather}
-\alpha\bigl(x^{\prime }(s)\sin \theta+y^{\prime }(s)\cos \theta\bigr)= \mathcal{A}\bigl(x^{\prime
}(s)\cos \theta-y^{\prime }(s)\sin \theta\bigr)+\mathcal{B}(-x(s)\sin \theta-y(s)\cos
\theta),\nonumber \\
\alpha\bigl(x^{\prime }(s)\cos \theta-y^{\prime }(s)\sin \theta\bigr)=\mathcal{A}\bigl(x^{\prime
}(s)\sin \theta+y^{\prime }(s)\cos \theta\bigr)+\mathcal{B}(x(s)\cos \theta-y(s)\sin \theta),\nonumber
\\
1=\mathcal{B}\bigl(x^2(s)+y^2(s)\bigr). \label{ab}
\end{gather}
The last equation of \eqref{ab} yields
$
\mathcal{B}=\frac{1}{x^2(s)+y^2(s)}$,
and hence
\[
b=\frac{\mathcal{B}}{\sqrt{1+\alpha^2}}=\frac{2\lambda^{2}}{\sqrt{1+\alpha^2}(1-\cos{2\lambda s})}.
\]
The first two equations of \eqref{ab} indicate
$
-\alpha\bigl(\bigl(x^{\prime 2}\bigr)+\bigl(y^{\prime 2}\bigr)\bigr)=\mathcal{B}\bigl(-xx^{\prime }-yy^{\prime }\bigr)$,
which implies
\begin{equation} \label{ac}
\alpha=\frac{xx^{\prime }+yy^{\prime }}{x^2+y^2}.
\end{equation}
In what follows, we claim the above $\alpha$ is one of special solutions.
Notice that \eqref{ad} shows
\[
x^2+y^2=\frac{1}{4\lambda^2}(2-2\cos(2\lambda s)),
\]
and hence \eqref{ac} can be rewritten as
\begin{equation*}
\alpha=\frac{1}{2}\bigl(\ln x^2+y^2\bigr)^{\prime }=\frac{\lambda \sin(2\lambda s)}{(1-\cos(2\lambda s))}.
\end{equation*}
Substituting $b$ and $\alpha$ into \eqref{ab}, we have
$
a=\frac{\mathcal{A}}{\sqrt{1+\alpha^2}}=\frac{-\lambda}{\sqrt{1+\alpha^2}}$.
\end{proof}

Given a $\alpha$-function, we have shown \cite{ChiuH/LiuH:2022} that the first fundamental form ($a$, $b$) is determined up to two functions $h(y)$ and $k(y)$ as follows.
\begin{Proposition}
\label{metricform} For any $\alpha(x,y)=\frac{c}{2}\frac{\sin{(cx+c_{1}(y)})
}{c_{2}(y)-\cos{(cx+c_{1}(y))}}$, the explicit formula for the induced
metric on a constant $p$-mean curvature surface with $c$ as its $p$-mean
curvature and this $\alpha$ as its $\alpha$-function is given by
\begin{equation*} \label{meta}
a=\left(-\frac{c}{2}+\frac{\frac{c}{2}c_{2}}{(c_{2}-\cos{(cx+c_{1})})}+\frac{
h(y)}{|c_{2}-\cos{(cx+c_{1})}|}\right)\frac{1}{\bigl(1+\alpha^{2}\bigr)^{1/2}},
\end{equation*}
and
\begin{equation*} \label{metb}
b=\left(\frac{e^{k(y)}}{|c_{2}-\cos{(cx+c_{1})}|}\right)\frac{1}{
\bigl(1+\alpha^{2}\bigr)^{1/2}},
\end{equation*}
for some functions $h(y)$ and $k(y)$.
\end{Proposition}

\subsection{The normalization} As we normalize the induced metric $a$ and $b$ to be close as much as possible to the metric induced on the horizontal $p$-minimal plane, we would like to normalize $a$ and $b$ so that they look like the induced
metric of the Pansu sphere. Indeed, from the transformation law \cite[formula~(2.20)]{ChiuH/LiuH:2022}, it is easy to see that there exist
another compatible coordinates $(\tilde{x},\tilde{y})$, called
normal coordinates such that
\begin{gather}
\tilde{a}=-\frac{\frac{c}{2}}{\bigl(1+\alpha^{2}\bigr)^{1/2}}, \label{nora}
\\
\tilde{b}=\frac{\frac{c^2}{2}}{|c_{2}-\cos{(cx+c_{1})}|\bigl(1+\alpha^{2}\bigr)^{1/2}}, \label{norb}
\end{gather}
where $c=2\lambda$. Such normal coordinates are uniquely determined up to a translation. We thus have the following theorem.

\begin{Theorem}
\label{cominv} In normal coordinates $(x,y)$, the functions $c_{1}(y)$ and $
c_{2}(y)$ in the expression \smash{$\alpha(x)=\frac{c}{2}\frac{\sin{(cx+c_{1}})}{
c_{2}-\cos{(cx+c_{1})}}$} are unique in the following sense: up to a
translation on $y$, $c_{2}(y)$ is unique, and $c_{1}(y)$ is unique up to a
constant. We denote these two unique functions by
$
\zeta_{1}(y)=c_{1}(y)$, $ \zeta_{2}(y)=c_{2}(y)$.
Therefore, the set $\{\zeta_{1}(y), \zeta_{2}(y)\}$ constitutes a \emph{complete} set of
invariants for those surfaces $(\alpha$ not vanishing$)$.
\end{Theorem}

It is worth our attention that, for the surfaces with $c_{2}>1$, the
denominator of the formula for $\alpha$ is never zero. That means the
surfaces won't extend to a surface with singular points. Moreover, if the
surface is closed, it must be a closed constant $p$-mean curvature
surface without singular points, which means the surface is of type of
torus. This indicates that it is possible to find a Wente-type torus in this
class of surfaces.

\subsection{The structure of the singular sets} In this subsection, we study the structure of the singular set. For the general type, we choose a~normal coordinate system $(s,\theta )$
such that
\begin{equation*}
\alpha =\lambda \frac{\sin \left( 2\lambda s+\zeta _{1}(\theta )\right) }{\zeta _{2}(\theta )-\cos \left( 2\lambda s+\zeta _{1}(\theta )\right) },
\end{equation*}
and
\begin{equation*}
a=-\frac{\lambda }{\sqrt{1+\alpha ^{2}}},\qquad b=\frac{2\lambda ^{2}}{\left \vert \zeta _{2}( \theta ) -\cos \left( 2\lambda s+\zeta
_{1}( \theta ) \right) \right \vert \sqrt{1+\alpha ^{2}}}.
\end{equation*}
Then the singular set is the graph of the function
\smash{$
x( \theta ) =\frac{\cos ^{-1}\left( \zeta _{2}(\theta )\right)-\zeta _{1}(\theta )}{2\lambda }$}.
The induced metric~$I$ (or the first fundamental form) on the regular part
reads
\begin{equation*}
I={\rm d}s\otimes {\rm d}s-\frac{a}{b}{\rm d}s\otimes {\rm d}\theta -\frac{a}{b}{\rm d}\theta \otimes {\rm d}s+
\frac{\bigl( 1+a^{2}\bigr) }{b^{2}}{\rm d}\theta \otimes {\rm d}\theta .
\end{equation*}
Now we use the metric to compute the length of the singular set \smash{$\big\{
\bigl( \frac{\cos ^{-1}( \zeta _{2}(\theta )) -\zeta _{1}(\theta )
}{2\lambda },\theta \bigr) \big\} $}, where~$\theta $ belongs to some
open interval.

Case $\zeta _{2}(\theta ) \neq 1$. Let \smash{$\gamma (\theta )=\bigl( \frac{\cos ^{-1}( \zeta _{2}(\theta ))
-\zeta _{1}(\theta )}{2\lambda },\theta \bigr) $}, which is a parametrization
of the singular set. Then the square of the velocity at $\theta $ is
\begin{align}
\left \vert \gamma ^{\prime }(\theta )\right \vert ^{2} =\left[ q^{\prime
}( \theta ) \right] ^{2}-\frac{2aq^{\prime }( \theta )
}{b}+\frac{a^{2}+1}{b^{2}} =\frac{\left[ a-bq^{\prime }( \theta ) \right] ^{2}+1}{b^{2}}>0
\qquad\text{for all}\quad\theta , \label{sglength}
\end{align}
where
\begin{equation*}
q^{\prime }( \theta ) = \left[ \frac{\cos ^{-1}( \zeta
_{2}(\theta )) -\zeta _{1}(\theta )}{2\lambda }\right] ^{\prime }
= \frac{-\zeta _{2}^{^{\prime }}( \theta ) }{2\lambda \sqrt{
1-\zeta _{2}^{2}( \theta ) }}-\frac{\zeta _{1}^{^{\prime }}(
\theta ) }{2\lambda }.
\end{equation*}
Formula \eqref{sglength} shows that the parametrized curve $\gamma (\theta )$
of the singular set has a positive length.

Case $\zeta _{2}(\theta ) =1$. We parametrize the singular set by \smash{$\gamma
(\theta )=\bigl( \frac{\cos ^{-1}( 1-\epsilon ) -\zeta
_{1}(\theta )}{2\lambda },\theta \bigr) $} for $\epsilon >0$. It is easy to see
\smash{$
\gamma ^{\prime }( \theta ) =\big( \frac{-\zeta _{1}^{\prime
}( \theta ) }{2\lambda },1\big) $}.
When $\varepsilon \rightarrow 0$, the metric
\[
I={\rm d}s\otimes {\rm d}s-\frac{a}{b}{\rm d}s\otimes {\rm d}\theta -\frac{a}{b}{\rm d}\theta \otimes {\rm d}s+
\frac{\bigl( 1+a^{2}\bigr) }{b^{2}}{\rm d}\theta \otimes {\rm d}\theta
\]
degenerates to
$\tilde{I}={\rm d}s\otimes {\rm d}s$.

Then the square of the velocity at $\theta $ is
\smash{$
\big|\gamma ^{\prime }( \theta ) \big|^{2}=\frac{[ \zeta _{1}^{\prime
}( \theta ) ] ^{2}}{4\lambda ^{2}}$}.
Thus, if $\zeta _{1}( \theta ) =c_{1}$ and $\zeta _{2}( \theta
) =1$, the length of the parametrized curve $\gamma (\theta )$ of
the singular set is zero. This result coincides with the singular set for the Pansu sphere being isolated.
We conclude the above discussion with the following theorem, an analog of \cite[Theorem 1.7]{ChiuH/LiuH:2022}.

\begin{Theorem}\label{strofsiset}
The singular set of a constant $p$-mean surface with $H=c\neq 0$ is either
\begin{enumerate}\itemsep=0pt
\item[$(1)$] an isolated point; or
\item[$(2)$] a smooth curve.
\end{enumerate}
In addition, an isolated singular point only happens on the surfaces of special type I with $\zeta_{1}=\text{\rm const}$, namely, a part of the Pansu sphere containing one of the poles as the isolated singular point.
\end{Theorem}
Theorem \ref{strofsiset} together with \cite[Theorem 1.7]{ChiuH/LiuH:2022} are just special cases of \cite[Theorem 3.3]{ChengJ/HwangJ/MalchiodiA/YangP:2005}. However, we give a computable proof of this result for constant $p$-mean surfaces. We also have the description of how a characteristic leaf goes through a singular curve, which is called a ``go through'' theorem in \cite{ChengJ/HwangJ/MalchiodiA/YangP:2005}. Suppose $p_{0}$ is a point in a singular curve. From the above basic facts, we see that a characteristic curve $\gamma$ always reaches the singular point $p_{0}$ going a finite distance. From the opposite direction, suppose $\tilde{\gamma}$ is another characteristic curve that reaches $p_{0}$. Then the union of $\gamma$, $p_{0}$ and $\tilde{\gamma}$ forms a smooth curve (we also refer the reader to the proof of \cite[Theorem~1.8]{ChiuH/LiuH:2022}, they are similar). We thus have the following theorem.

\begin{Theorem}\label{theo2}
Let $\Sigma\subset H_{1}$ be a constant $p$-mean surface with $H=c\neq 0$. Then the characteristic foliation is smooth around the singular curve in the following sense that each leaf can be extended smoothly to a point on the singular curve.
\end{Theorem}

Making use of Theorem \ref{theo2}, we have the following result.
\begin{Theorem}\label{main10}
Let $\Sigma$ be a constant $p$-mean surface of type II $($III$)$ with $H=c\neq 0$. If it can be smoothly extended through the singular curve, then the other side of the singular curve is of type III $($II$)$.
\end{Theorem}

Therefore, we see that a surface of general type $II$ is always pasted together with a surface of general type $III$ at a singular curve and vice versa.

\section[Rotationally invariant surfaces in H\_1]{Rotationally invariant surfaces in $\boldsymbol{ H_1}$}

Let $\Sigma(s,\theta)$ be a rotationally invariant surface in $H_1$
generated by a curve $\gamma(s)=(x(s),0,t(s))$ on the $xt$-plane,
that is, $\Sigma$ is parametrized by
$
\Sigma(s,\theta)=(x(s)\cos \theta,x(s)\sin \theta,t(s))$,
where ${x'^{2}+t'^{2}=1}$. Here $'$ means taking a derivative with respect to $s$.

\subsection[The computation of H, alpha, a and b]{The computation of $\boldsymbol{H}$, $\boldsymbol{\alpha}$, $\boldsymbol{a}$ and $\boldsymbol{b}$}

Now we consider the horizontal (see \cite[Definition 1.1]{ChiuH/HuangY/LaiS:2017}) generating curve
\[
\tilde{\gamma}(s)=(x(s)\cos \theta(s),x(s)\sin \theta(s),t(s)).
\]

\begin{Lemma}
$\tilde{\gamma}$ is horizontal if and only if $t^{\prime}+x^{2}\theta^{\prime}=0$.
\end{Lemma}

\begin{proof}
Note that at the point $\tilde{\gamma}(s)$,
\begin{gather*}
\mathring e_1=\frac{\partial}{\partial x_1}+y_1\frac{\partial}{\partial z}=
\frac{\partial}{\partial x_1}+x(s)\sin \theta(s)\frac{\partial}{\partial z},
\\
\mathring e_2=\frac{\partial}{\partial y_1}-x_1\frac{\partial}{\partial z}=
\frac{\partial}{\partial y_1}-x(s)\cos \theta(s)\frac{\partial}{\partial z},
\end{gather*}
and direct computations imply
\begin{equation*}
\tilde{\gamma}^{\prime }(s)=\big(x^{\prime }\cos \theta-x\theta^{\prime }\sin
\theta\big)\mathring e_1+\big(x^{\prime }\sin \theta+x\theta^{\prime }\cos
\theta\big)\mathring e_2+\big(t^{\prime}+x^{2}\theta^{\prime}\big)T,
\end{equation*}
and hence $\tilde{\gamma}^{\prime }(s)\in \xi$ if and only if $t^{\prime}+x^{2}\theta^{\prime}=0$.
\end{proof}

Let $\tilde{s}$ be the horizontal arc-length of $\tilde{\gamma}(s)$. We can thus re-parametrize the surface $\Sigma(s,\theta)$ to~be
\begin{align*} 
\Sigma(\tilde s,\tilde \theta)={}&\bigl(x(s)\cos \theta(s)\cos \tilde
\theta-x(s)\sin \theta(s)\sin \tilde \theta,x(s)\cos \theta(s)\sin \tilde
\theta\\
&+x(s)\sin \theta(s)\cos \tilde \theta,t(s)\bigr),
\end{align*}
with a compatible coordinate system
\[
e_1=\frac{\partial}{\partial \tilde s}=\Sigma_s\frac{\partial s}{\partial
\tilde s}, \qquad \text{where} \quad \Sigma_s=\bigl(x^{\prime }\cos
\phi-x\theta^{\prime }\sin \phi\bigr)\mathring e_1+
\bigl(x^{\prime }\sin \phi+x\theta^{\prime }\cos \phi\bigr)\mathring e_2.
\]
Moreover, we see
\smash{$
|\tilde{\gamma}|^{\prime 2}=\frac{x^2x^{\prime 2}+t^{\prime 2}}{x^2}$},
so that we may choose $\tilde s$ such that \smash{$\big|\frac{{\rm d}\tilde{\gamma}({\tilde s}%
)}{d\tilde{s}}\big|=1$}, that is,
\begin{equation}\label{excss}
\frac{{\rm d}\tilde s}{{\rm d}s}=|\tilde{\gamma}^{\prime }(s)|=\frac{\sqrt{x^2x^{\prime 2}+t^{\prime 2}}}{x}.
\end{equation}

Manipulating $\Sigma$ to be
\begin{equation*} \label{surface2}
\Sigma\bigl(\tilde s,\tilde \theta\bigr)=\bigl(x(s)\cos \bigl(\theta(s)+\tilde \theta\bigr),
x(s)\sin\bigl( \theta(s)+\tilde \theta\bigr), t(s)\bigr),\qquad \bigl(\mbox{denote }
\phi=\theta(s)+\tilde \theta\bigr)
\end{equation*}
and obtain
\begin{gather*}
\Sigma_{\tilde s}=\frac{{\rm d} s}{{\rm d} \tilde s}\bigl(x^{\prime }\cos
\phi-x\theta^{\prime }\sin \phi\bigr)\mathring e_1+\frac{{\rm d}s}{{\rm d} \tilde s}
\bigl(x^{\prime }\sin \phi+x\theta^{\prime }\cos \phi\bigr)\mathring e_2, \\
\Sigma_{\tilde \theta}=-x\sin \phi \mathring e_1+x\cos \phi \mathring
e_2+x^2\frac{\partial}{\partial z}.
\end{gather*}
Then
\begin{gather}
e_1=\Sigma_{\tilde s}=\frac{x}{\sqrt{x^2x^{\prime 2}+t^{\prime 2}}}
\bigl(x^{\prime }\cos \phi-x\theta^{\prime }\sin \phi\bigr)\mathring e_1+\frac{{\rm d}s}{{\rm d}
\tilde s}\bigl(x^{\prime }\sin \phi+x\theta^{\prime }\cos \phi\bigr)\mathring e_2, \label{e1} \\
e_2=Je_1=\frac{x}{\sqrt{x^2x^{\prime 2}+t^{\prime 2}}}\bigl(x^{\prime }\cos
\phi-x\theta^{\prime }\sin \phi\bigr)\mathring e_2-\frac{{\rm d}s}{{\rm d} \tilde s}
\bigl(x^{\prime }\sin \phi+x\theta^{\prime }\cos \phi\bigr)\mathring e_1. \label{e2}
\end{gather}
The fact that \smash{$\frac{\alpha e_2+T}{\sqrt{1+\alpha^2}} \in T\Sigma$} implies
$
\alpha e_2+T=a\sqrt{1+\alpha^2}\Sigma_{\tilde s}+b\sqrt{1+\alpha^2}
\Sigma_{\tilde \theta}$.
Using \eqref{e1}, \eqref{e2} and comparing the coefficients of $\mathring
e_1$, $\mathring e_2$, and $T$, respectively, one sees
\begin{equation} \label{a7}
a=\frac{t^{\prime }}{x\sqrt{1+\alpha^2}\sqrt{x^2x^{\prime 2}+t^{\prime 2}}}
,\qquad b=\frac{1}{x^2\sqrt{1+\alpha^2}},\qquad \alpha=\frac{x^{\prime }}{
\sqrt{x^2x^{\prime 2}+t^{\prime 2}}},
\end{equation}
and hence, from the first equation of the integrability conditions \eqref{intcon}, we have
\begin{equation} \label{a8}
H=-\frac{x^{3}\bigl(x^{\prime }t^{\prime \prime }-x^{\prime \prime
}t^{\prime }\bigr)+t^{\prime 3}}{x\big\{x^2x^{\prime 2}+t^{\prime 2}\big\}^{3/2}}.
\end{equation}

\subsection[Another understanding of energy E]{Another understanding of energy $\boldsymbol{ E}$}

In this subsection, we assume moreover that the rotationally invariant surface $\Sigma$ is of constant $p$-mean curvature. We consider the relation between the integrability condition and the
energy discussed in Ritor\'e and Rosales' paper \cite{Ritore/Rosales:2006}. The integrability condition
$
a_{\tilde s}-\frac{a}{b}b_{\tilde s}+\frac{c\alpha}{\sqrt{1+\alpha^2}}=0
$
indicates that
\begin{equation} \label{a5}
\int \left(\frac{a_{\tilde s}}{b}-\frac{a}{b^2}b_{\tilde s}+\frac{c\alpha}{b\sqrt{
1+\alpha^2}}\right){\rm d}\tilde s
\end{equation}
is a constant. Then we have \eqref{a5} computed as
\begin{equation*}
\frac{a}{b}+c\int \frac{x^2x^{\prime }}{\sqrt{x^2x^{\prime 2}+t^{\prime 2}}}
{\rm d}\tilde s
=\frac{a}{b}+c\int xx^{\prime }{\rm d}s
=\frac{a}{b}+\lambda x^2, \qquad\textrm{up to a constant},
\end{equation*}
which clearly says that $\frac{a}{b}+\lambda x^2$ is constant.
The constant $\frac{a}{b}+\lambda x^2$ interprets the energy $E$ based on
Ritor\'e's discussion. Indeed, we have
\begin{equation} \label{ffe}
E=\frac{a}{b}+\lambda x^2 =\frac{xt^{\prime }}{\sqrt{x^2x^{\prime 2}+t^{\prime 2}}}+\lambda x^2 =t_{\tilde s}+\lambda x^2,
\end{equation}
that is, $t_{\tilde s}=E-\lambda x^2$. One sees
\begin{equation} \label{function_t}
t=E\tilde s-\lambda \int x^2{\rm d}\tilde s.
\end{equation}

\subsection{The Coddazi-like equation}

For later use, we calculate $1+\alpha^2=\frac{1+x^2x^{\prime 2}}{x^2x^{\prime 2}+t^{\prime 2}}$ and convert $\alpha$ to be of the general form
\begin{equation} \label{a1}
\alpha=\frac{x^{\prime }}{\sqrt{x^2x^{\prime 2}+t^{\prime 2}}}=\frac{
x_{\tilde s}}{x}.
\end{equation}
Note that $\alpha$ satisfies the Coddazi-like equation \smash{$\alpha_{\tilde s\tilde s}+6\alpha\alpha_{\tilde s}+4 \alpha^{3}+c^2\alpha=0$}, where $c=2\lambda$. Then this ODE immediately shows
\begin{equation} \label{a2}
\frac{x_{\tilde s\tilde s\tilde s}}{x}+3\frac{x_{\tilde s}x_{\tilde s\tilde
s}}{x^2}+c^2\frac{x_{\tilde s}}{x}=0.
\end{equation}
The equation \eqref{a2} is manipulated to be
$
\bigl(xx_{\tilde s\tilde s}+(x_{\tilde s})^2+\frac{c^2}{2}x^2\bigr)_{\tilde s}=0$,
which gives
\begin{equation} \label{a3}
\big(x^2\big)_{\tilde s\tilde s}+c^2x^2=k,\qquad \mbox{for some constant $k$.}
\end{equation}
Let $u=x^2$, then \eqref{a3} becomes a second-order inhomogeneous constant
coefficient ODE
\begin{equation} \label{a4}
u_{\tilde s\tilde s}+c^2u=k.
\end{equation}
\begin{itemize}\itemsep=0pt
\item[(I)] Suppose $c\neq 0$, the homogeneous ODE $u_{\tilde s\tilde s}+c^2u=0$ has the
general solution $u_h$ given by
\begin{equation}\label{rk1}
u_h=k_1\sin(c\tilde s)+k_2\cos(c\tilde s)=r\cos(c\tilde s-c_1),
\end{equation}
where $r=\sqrt{k_1^2+k_2^2}$ and $r\sin c_1=k_1$. One also notes that $u_p=%
\frac{k}{c^2}$ is a particular solution to~\eqref{a4}, and hence
\begin{equation} \label{solution_u}
x^2=u=\frac{k}{c^2}+r\cos(c\tilde s-c_1).
\end{equation}

\item[(II)] When $c=0$, it is clear that \eqref{a3} becomes $\big(x^2\big)_{\tilde s\tilde s}=k$, which implies
\begin{equation}\label{solution_x}
x^2=k\tilde s^2+2k_1\tilde s+k_2,
\end{equation}
for some constants $k$, $k_{1}$ and $k_{2}$.
\end{itemize}
\begin{Example}
If $k=0$, $k_{1}=0$, then \eqref{solution_x} yields $x=\sqrt{k_2}$. On
the other hand, \eqref{pmce} suggests~\smash{$
0=\frac{t^{\prime 3}}{\sqrt{k_2} t^{\prime 3}}=\frac{1}{\sqrt{k_2}}>0$},
which is a contradiction. We conclude that there are no such kinds of $p$-minimal surfaces ($k=0$, $k_{1}=0$) which are rotationally symmetric. In this case, $\alpha$~vanishes so that it corresponds a~vertical cylinder surface which is absolutely not $p$-minimal.
\end{Example}

\subsection[The relation between k and E]{The relation between $\boldsymbol{k}$ and $\boldsymbol{E}$}
Assume that $c=2\lambda\neq 0$. We write \eqref{a8} as
\begin{align}
-2\lambda&=\frac{x^{3}\bigl(x^{\prime }t^{\prime \prime }-x^{\prime \prime
}t^{\prime }\bigr)+t^{\prime 3}}{x\big\{x^2x^{\prime 2}+t^{\prime 2}\big\}^{3/2}}\nonumber\\
&=\frac{x^{2}\bigl(\frac{t'}{x'}\bigr)'x^{\prime 2}}{\big\{x^2x^{\prime 2}+t^{\prime 2}\big\}^{3/2}}+\frac{1}{x^4}\left(\frac{xt^{\prime}}{\sqrt{x^{2}x^{\prime 2}+t^{\prime 2}}}\right)^3=I_1+I_2,\label{pmce}
\end{align}
where \[I_1=\frac{x^{2}\bigl(\frac{t'}{x'}\bigr)'x^{\prime 2}}{\big\{x^2x^{\prime 2}+t^{\prime 2}\big\}^{3/2}},\qquad I_2=\frac{1}{x^4}\left(\frac{xt^{\prime}}{\sqrt{x^{2}x^{\prime 2}+t^{\prime 2}}}\right)^3.\]
From \eqref{solution_u}, taking a derivative with respect to $s$, we have
\begin{equation}\label{solution_ud}
x'=\frac{-cr\sin{(c\tilde{s}-c_1)}\sqrt{x^2x^{\prime 2}+t^{\prime 2}}}{2x^2}.
\end{equation}
On the other hand, \eqref{ffe} implies
\begin{equation}\label{ffed}
t'=\bigl(E-\lambda x^2\bigr)\frac{\sqrt{x^2x^{\prime 2}+t^{\prime 2}}}{x}.
\end{equation}
By means of \eqref{solution_u}, \eqref{solution_ud},\eqref{ffed} and \eqref{excss}, after direct computations, we have
\begin{align*}
I_1={}&\frac{x^{2}\bigl(\frac{t'}{x'}\bigr)'x^{\prime 2}}{\big\{x^2x^{\prime 2}+t^{\prime 2}\big\}^{3/2}}\\
={}&\left(\frac{2x\bigl(E-\lambda x^2\bigr)}{-cr\sin{(c\tilde{s}-c_1)}}\right)^{\prime}\left(\frac{c^{2}r^{2}\sin^{2}{(c\tilde{s}-c_1)}\bigl(x^2x^{\prime 2}+t^{\prime 2}\bigr)}{4x^4}\right)\\
={}&\frac{1}{4x^2}\Biggl(-2x^{2}x'\bigl(4\lambda x^{2}x'-2x'\bigl(E-\lambda x^{2}\bigr)\bigr)\left(\frac{x^{2}-\bigl(E-\lambda x^{2}\bigr)^{2}}{x^{4}x^{\prime 2}}\right) \\
& +2\bigl(E-\lambda x^{2}\bigr)\bigl(c^{2}x^{2}-k\bigr)\Biggr),
\end{align*}
and from \eqref{ffe}, we have
\begin{equation*}
I_2=\frac{1}{x^4}\left(\frac{xt^{\prime}}{\sqrt{x^{2}x^{\prime 2}+t^{\prime 2}}}\right)^3=\frac{1}{x^{4}}\bigl(E-\lambda x^{2}\bigr)^3.
\end{equation*}
Therefore,
\[
-2\lambda=I_1+I_2=\frac{1}{4x^4}\bigl(\bigl(-4c-2Ec^{2}-2c+ck\bigr)x^4+\bigl(8\lambda E^{2}+4E-2Ek\bigr)x^2\bigr),
\]
which implies
$
k=2cE+2$.

\subsection[Horizontal generating curves for c neq 0]{Horizontal generating curves for $\boldsymbol{c\neq 0}$}

In this subsection, we will show that $k$, and hence the energy $E$, has a lower bound. A horizontal generating curve of a rotationally invariant constant $p$-mean curvature surface is a geodesic curve, which is parametrized by
\begin{equation}\label{gd}
\tilde{\gamma}(\tilde{s})=\left(\begin{matrix}
\dfrac{1}{c}\sin(c\tilde{s})+x_{0}\vspace{1mm}\\
-\dfrac{1}{c}\cos(c\tilde{s})+\frac{1}{c}+y_{0} \vspace{1mm}\\
\left(\dfrac{1}{c^2}+\dfrac{cy_{0}}{c^2}\right)\sin(c\tilde{s})+\dfrac{x_{0}}{c}\cos(c\tilde{s})-\dfrac{\tilde{s}}{c}+\dfrac{\pi}{c^2}-\dfrac{x_{0}}{c}+t_{0}\end{matrix} \right)
\end{equation}
for some $(x_{0},y_{0},t_{0})$, where $\tilde{s}$ is a horizontal arc length parameter.

Suppose that $\gamma(s)=(x(s),0,t(s))$ with $x\geq 0$ is the corresponding generating curve, we have
\begin{align*}
x^{2}&=\left(\frac{1}{c}\sin(c\tilde{s})+x_{0}\right)^{2}+\left(-\frac{1}{c}\cos(c\tilde{s})+\frac{1}{c}+y_{0}\right)^{2}\\
&=2\frac{x_{0}}{c}\sin(c\tilde{s})-2\left(\frac{1+cy_{0}}{c^2}\right)\cos(c\tilde{s})+x_{0}^{2}+\left(\frac{1+cy_{0}}{c}\right)^{2}+\frac{1}{c^2}\\
&=r\cos(c\tilde{s}-c_{1})+\frac{k}{c^2},
\end{align*}
 where
 \begin{align*}
 k&=1+(cx_{0})^{2}+(1+cy_{0})^{2}\geq 1,\\
 r&=\sqrt{\left(2\frac{x_{0}}{c}\right)^{2}+\left(2\left(\frac{1+cy_{0}}{c^2}\right)\right)^{2}}=\frac{2}{c^2}\sqrt{(cx_{0})^{2}+(1+cy_{0})^{2}}=\frac{2}{c^2}\sqrt{k-1},
 \end{align*}
 and $c_{1}$ is a real number such that
\smash{$\sin{c_{1}}=\frac{2\frac{x_{0}}{c}}{r}$}, \smash{$ \cos{c_{1}}=-\frac{2(\frac{1+cy_{0}}{c^2})}{r}$}.

\subsection[The invariants zeta\_1 and zeta\_2 for surfaces with c neq 0]{The invariants $\boldsymbol{\zeta_{1}}$ and $\boldsymbol{\zeta_{2}}$ for surfaces with $\boldsymbol{c\neq 0}$}
If $r=0$, then $k=1$, $x_0=0$, $y_0=-\frac{1}{c}$. Thus, \eqref{gd} implies
 \[
 \tilde\gamma(\tilde s)=\left(\frac{1}{c}\sin(c\tilde s),-\frac{1}{c}\cos(c\tilde s),-\frac{\tilde s}{c}+\frac{\pi}{c^2}+t_0\right),
 \]
 which generates a cylinder. We assume from now on that $r\neq 0$. Taking the derivative with respect to $\tilde s$ on both sides of
\eqref{solution_u} to have
$
2xx_{\tilde s}=-rc\sin(c\tilde s-c_1)$.
Together with \eqref{a1}, we have $\alpha$ of the general form as follows:
\begin{equation*}
\alpha=\frac{xx_{\tilde s}}{x^2}=\frac{-r\lambda \sin(c\tilde s-c_1)}{\frac{k}{c^2}+r\cos(c\tilde s-c_1)}=\frac{\lambda \sin(c\tilde s-c_1)}{c_2-\cos(c\tilde s-c_1)},
\end{equation*}
where $c_2=-\frac{k}{rc^2}$.

In this subsection, we want to normalize $a$ and $b$ such that they have the
forms looking as~\eqref{nora} and \eqref{norb}, respectively. Together with %
\eqref{a7}, \eqref{solution_u} and \eqref{ffe}, we have
\begin{align*}
a&=\frac{t^{\prime }}{x\sqrt{1+\alpha^2}\sqrt{x^2x^{\prime 2}+t^{\prime 2}}}=\left(\frac{E}{x^{2}}-\lambda \right)\frac{1}{\sqrt{1+\alpha^{2}}}, \\
b&=\frac{1}{x^2\sqrt{1+\alpha^2}}=\frac{1}{\frac{k}{c^2}+r\cos(c\tilde s-c_1)%
}\frac{1}{\sqrt{1+\alpha^{2}}}.
\end{align*}

Thus we choose the normal coordinates $\big\{ \bar{s},\bar{\theta}\big\}$ with
$
\bar{s}=\tilde{s}+\Gamma{\bigl(\tilde{\theta}\bigr)}$, $ \bar{\theta}=\Psi{\bigl(\tilde{%
\theta}\bigr)}$,
such that
$
\Gamma^{\prime }{\bigl(\tilde{\theta}\bigr)}=-E$, $ \Psi^{\prime }{\bigl(\tilde{\theta}\bigr)}%
=-2\lambda^{2}r$.
Then we have
\begin{equation*}
\bar{a}=\frac{-\lambda}{\sqrt{1+\bar \alpha^{2}}}, \qquad
\bar{b}=\frac{2\lambda^{2}}{\bigl(-\frac{k}{rc^2}-\cos(c\tilde s-c_1)\bigr)%
}\frac{1}{\sqrt{1+\bar \alpha^{2}}},
\end{equation*}
with
\[
\bar \alpha=\frac{\lambda \sin(c\tilde s-c_1)}{c_2-\cos(c\tilde s-c_1)}=%
\frac{\lambda \sin \bigl(c\bar s-c_1-\frac{E\bar \theta}{r\lambda}\bigr)}{%
c_2-\cos \bigl(c\bar s-c_1-\frac{E\bar \theta}{r\lambda}\bigr)},
\]
that is,
\begin{equation}\label{coinro}
\zeta_{1}\bigl(\bar \theta\bigr)=-c_1-\frac{2E\bar \theta}{cr},\qquad
\zeta_{2}\bigl(\bar \theta\bigr)=c_{2}=-\frac{k}{rc^2}=-\frac{2cE+2}{c^{2}r}.
\end{equation}

If $E=0$, then $k=2$, thus the surface has the generating curve defined by
\begin{gather*}
x^{2}=\frac{2}{c^2}+r\cos{(c\tilde{s}-c_1)}=r\left(\frac{2}{c^{2}r}+\cos{(c\tilde{s}-c_1)}\right),\\
t=-\lambda\left(\frac{2}{c^2}\tilde{s}+\frac{r}{c}\sin{(c\tilde{s}-c_1)}\right),
\end{gather*}
with
$\zeta_{1}\bigl(\bar \theta\bigr)=-c_{1}$, $ \zeta_{2}\bigl(\bar \theta\bigr)=-\frac{2}{c^{2}r}<0$.
Therefore, we see that
$x^{2}\geq 0\Leftrightarrow \cos{(c\tilde{s}-c_1)}\geq\zeta_{2}\bigl(\bar \theta\bigr)$,
which means that the generating curve $(x,t)$ is defined on the whole $\R$ if and only if $\zeta_{2}\bigl(\bar \theta\bigr)\leq-1$.
In particular, if $\zeta_{2}\bigl(\bar \theta\bigr)=-1$, it is the Pansu sphere.

If $E\neq 0$, then $k=2cE+2$ and
 \eqref{coinro} implies that
\begin{equation} \label{coinro1}
\zeta_{1}'\bigl(\bar \theta\bigr)=-\frac{2E}{cr}>\zeta_{2}\bigl(\bar \theta\bigr).
\end{equation}
For any constants $\eta_{1}$ and $\eta_{2}$ with $\eta_{1}>\eta_{2}$, we
obtain the unique solution to the equation system
\[
-\frac{2E}{c r}=\eta_{1},\qquad -\frac{2cE+2}{c^{2}r}=\eta_{2}.
\]

\subsection[The allowed values of k and E with c=0]{The allowed values of $\boldsymbol{ k}$ and $\boldsymbol{E}$ with $\boldsymbol{c=0}$}

In this subsection, we shall show what possible values can $k$ and $E$ attain.
Assume that $c=2\lambda=0$. We write \eqref{a8} as
\begin{equation*}
0=\frac{x^{3}\bigl(x^{\prime }t^{\prime \prime }-x^{\prime \prime
}t^{\prime }\bigr)+t^{\prime 3}}{x\big\{x^2x^{\prime 2}+t^{\prime 2}\big\}^{3/2}}=\frac{x^{2}\bigl(\frac{t'}{x'}\bigr)'x^{\prime 2}}{\big\{x^2x^{\prime 2}+t^{\prime 2}\big\}^{3/2}}+\frac{1}{x^4}\left(\frac{xt^{\prime}}{\sqrt{x^{2}x^{\prime 2}+t^{\prime 2}}}\right)^3=I_1+I_2.
\end{equation*}
Taking a derivative of \eqref{solution_x} with respect to $s$, we get
\begin{equation}\label{solution_xd}
x'=\frac{(k\tilde{s}+k_{1})\sqrt{x^2x^{\prime 2}+t^{\prime 2}}}{x^2}.
\end{equation}
On the other hand, from \eqref{ffe}, we have
\begin{equation}\label{ffexd}
t'=\frac{E\sqrt{x^2x^{\prime 2}+t^{\prime 2}}}{x}.
\end{equation}
By means of \eqref{solution_x}, \eqref{solution_xd},\eqref{ffexd} and \eqref{excss}, a direct computation gives
\[
I_1=\frac{x^{2}\bigl(\frac{t'}{x'}\bigr)'x^{\prime 2}}{\big\{x^2x^{\prime 2}+t^{\prime 2}\big\}^{3/2}}=\frac{1}{x^2}\left(\frac{E\bigl(x^{2}-E^{2}\bigr)}{x^2}-kE\right),
\]
and \eqref{ffe} implies
\begin{equation*}
I_2=\frac{1}{x^4}\left(\frac{xt^{\prime}}{\sqrt{x^{2}x^{\prime 2}+t^{\prime 2}}}\right)^3=\frac{1}{x^{4}}E^3.
\end{equation*}
Therefore,
$0=I_1+I_2=\frac{E(1-k)}{x^2}$,
which says that
\begin{equation}\label{ekvl}
E=0\qquad \textrm{or}\qquad k=1.
\end{equation}
The equation \eqref{ffe} says that $t $ and $E\tilde s$ are differed only by a constant. If $E=0$, then $t$ is constant, which gives us that $\Sigma$ is a plane that is perpendicular to $t$-axis.

\subsection[The invariants zeta\_1 and zeta\_2 for surfaces with c=0 and E neq 0]{The invariants $\boldsymbol{\zeta_{1}}$ and $\boldsymbol{\zeta_{2}}$ for surfaces with $\boldsymbol{c=0}$ and $\boldsymbol{E\neq 0}$}
 If $E=0$, the surface is a perpendicular plane to the $t$-axis. Therefore, in this subsection, we assume that $E\neq 0$, and thus, from \eqref{ekvl}, we have $k=1$.
Then one rewrites $\alpha$ to be
\begin{equation*}
\alpha=\frac{xx_{\tilde s}}{x^2}=\frac{\tilde s+k_1}{(\tilde s+k_1)^2+\bigl(k_2-k_{1}^{2}\bigr)}.
\end{equation*}
From \eqref{solution_x}, we have
$
\frac{1}{x^2}=\frac{\alpha}{\tilde{s}+k_1}>0$.
We want to normalize $a$ and $b$ such that they have the form specified in \cite[Theorem~1.3]{ChiuH/LiuH:2022}.
Together with \eqref{a7} and \eqref{ffe}, we have
\begin{align*}
a&=\frac{t^{\prime }}{x\sqrt{1+\alpha^2}\sqrt{x^{2}x^{\prime 2}+t^{\prime 2}}}
=\left(\frac{E}{x^{2}}\right)\frac{1}{\sqrt{1+\alpha^{2}}}=Eb, \\
b&=\frac{1}{x^2\sqrt{1+\alpha^2}}=\frac{\alpha}{(\tilde{s}+k_1)\sqrt{1+\alpha^2}}=\frac{|\alpha|}{(|\tilde{s}+k_1|)\sqrt{1+\alpha^2}}.
\end{align*}
Thus we choose the normal coordinates $\big\{ \bar{s},\bar{\theta}\big\}$ with
$
\bar{s}=\tilde{s}+\Gamma{\bigl(\tilde{\theta}\bigr)}$, $ \bar{\theta}=\Psi{\bigl(\tilde{\theta}\bigr)}
$
such that
$
\Gamma^{\prime }{\bigl(\tilde{\theta}\bigr)}=-E$, $ \Psi^{\prime }{\bigl(\tilde{\theta}\bigr)}=1$.
Then we have
\begin{equation*}
\bar{a}=0,\qquad \bar{b}=\frac{|\bar \alpha|}{|\tilde{s}+k_1|\sqrt{1+\bar
\alpha^2}}
\end{equation*}
with
\[
\bar \alpha=\frac{\tilde s+k_1}{(\tilde s+k_1)^2+\bigl(k_2-k_{1}^{2}\bigr)}=\frac{\bar{s}-\Gamma{\bigl(
\tilde{\theta}\bigr)}+k_1}{(\bar{s}-\Gamma{\bigl(\tilde{\theta}\bigr)}+k_1)^2+\bigl(k_2-k_{1}^{2}\bigr)},
\]
that is,
\begin{alignat*}{3}
&\zeta_{1}\bigl(\bar \theta\bigr)=k_1+E\bar\theta,\qquad && \textrm{which is linear in} \ \bar{\theta},&\\
&\zeta_{2}\bigl(\bar \theta\bigr)=k_2-k_{1}^{2},\qquad && \textrm{which is a constant, denoted as} \ \zeta_{2}.&
\end{alignat*}

From \eqref{solution_x} and \eqref{ffe}, we conclude that the generating curve is defined by
\begin{equation}\label{gc}
x^{2}=(\tilde{s}+k_1)^{2}+\zeta_{2},\qquad
t=E\tilde{s},\qquad \textrm{up to a constant}.
\end{equation}

\begin{Remark} We remark that for $\lambda=0$, two kinds of $p$-minimal surfaces are presented depending on the energy $E$. When $E=0$, $t$ in (3.16) is constant and then one obtains a plane that is perpendicular to the $t$-axis. On the other hand, if $E\neq 0$, we have $p$-minimal surfaces generated by curves defined by \eqref{gc}. For $\lambda \neq 0$, substituting \eqref{solution_u} in \eqref{function_t}, we
see
\begin{align*}
t&=E\tilde s-\lambda \int \left(\frac{k}{c^2}+\sqrt{k_1^2+k_2^2}\cos(c\tilde
s-c_1)\right){\rm d}\tilde s \\
&=\left(E-\frac{k}{4\lambda}\right)\tilde s-\frac{\sqrt{k_1^2+k_2^2}}{2}\sin(2\lambda
\tilde s-c_1)+\text{\rm const}.
\end{align*}
\end{Remark}
In the case $\lambda \neq 0$, we give the following two examples.

\begin{Example}
We choose $k_1$, $k_2$ in \eqref{rk1} so that $\sqrt{k_1^2+k_2^2}=-\frac{k}{c^2}$, and then
\eqref{solution_u} implies
\begin{equation*}
x=\pm \frac{\sqrt{2k}}{2\lambda} \sin\left(\lambda \tilde s-\frac{c_1}{2}\right).
\end{equation*}
Moreover, if $E=0$, then $t=-\frac{k}{4\lambda}\tilde s+\frac{k}{8\lambda^2}
\sin(2\lambda \tilde s-c_1)$, which is a scaling sphere.
\end{Example}

The other two integrability conditions (see \cite[equation~(2.13)]{ChiuH/LiuH:2022}) are
\begin{equation} \label{a6}
-\frac{b_{\tilde s}}{b} = 2\alpha+\frac{\alpha \alpha_{\tilde s}}{1+\alpha^2
}, \qquad aH_{\tilde s}+bH_{\tilde \theta} =\frac{\alpha_{\tilde s\tilde s}+6\alpha
\alpha_{\tilde s}+4\alpha^3+\alpha H^2}{\sqrt{1+\alpha^2}}.
\end{equation}
We rewrite the first equation in \eqref{a6} as
\begin{equation*}
2\alpha+\frac{\alpha \alpha_{\tilde s}}{1+\alpha^2}+\frac{b_{\tilde s}}{b}=0.
\end{equation*}
Integrating on both sides to see that
\[
\int \left(2\alpha+\frac{\alpha \alpha_{\tilde s}}{1+\alpha^2}+\frac{b_{\tilde s}}{
b}\right){\rm d} \tilde s
\]
is a constant. More precisely, in terms of $x$, $x^{\prime }$, $t$, $t^{\prime }$, we
write
\begin{equation*}
\int 2\left(\alpha+\frac{\alpha \alpha_{\tilde s}}{1+\alpha^2}+\frac{b_{\tilde s}
}{b}\right){\rm d}\tilde s =\int \left(2\frac{x_{\tilde s}}{x}+\frac{\alpha \alpha_{\tilde s}}{1+\alpha^2}+
\frac{b_{\tilde s}}{b}\right){\rm d}\tilde s =\ln\bigl(bx^2\sqrt{1+\alpha^2}\bigr)+\text{\rm const}.
\end{equation*}
The conclusion is that $\ln\bigl(bx^2\sqrt{1+\alpha^2}\bigr)$ is a constant, which
also follows from \eqref{a7}.

Suppose $H$ is constant. The second equation of \eqref{a6} is exactly $
\alpha_{\tilde s\tilde s}+6\alpha \alpha_{\tilde s}+4\alpha^3+\alpha H^2=0$.
Using \eqref{a1}, this ODE becomes \eqref{a2}, which has been discussed
previously.

\section[The construction of constant p-mean curvature surfaces]{The construction of constant $\boldsymbol{ p}$-mean curvature surfaces}\label{ccpms}

In this section, we construct constant $p$-mean curvature
surfaces by perturbing the Pansu sphere in some way. Recall the
parametrization of the Pansu sphere \eqref{paps}. For each fixed angle $
\theta$, the curve $l_{\theta}$ defined by $l_{\theta}(s)=(x(s)\cos{\theta}
-y(s)\sin{\theta},x(s)\sin{\theta}+y(s)\cos{\theta},t(s))$ is a geodesic
with curvature $2\lambda$. Let $\mathcal{C}$ be an arbitrary curve $\mathcal{
C}\colon\R \rightarrow H_{1}$ given by $\mathcal{C}(\theta)=(x_{1}(\theta),x_{2}(
\theta),x_{3}(\theta))$. For each fixed $\theta$, we translate $l_{\theta}$
by $\mathcal{C}(\theta)$, so that the curve $L_{\mathcal{C}
(\theta)}(l_{\theta})$ is also a geodesic curve with curvature $2\lambda$.
Then the union of all these curves
$
\Sigma_{\mathcal{C}}=\cup_{\theta}L_{\mathcal{C}(\theta)}(l_{\theta})
$
constitutes \emph{a constant $p$-mean curvature surface} with a
parametrization
\begin{align}
Y(s,\theta)={}&(x_{1}(\theta)+(x(s)\cos{\theta}-y(s)\sin{\theta}\bigr)
,x_{2}(\theta)+(x(s)\sin{\theta}+y(s)\cos{\theta}),\nonumber \\
&x_{3}(\theta)+t(s)+x_{2}(\theta)(x(s)\cos{\theta}-y(s)\sin{\theta})\nonumber\\
&-x_{1}(\theta)(x(s)\sin{\theta}+y(s)\cos{\theta})).\label{pacms}
\end{align}
By a straightforward computation, and notice that
\begin{gather*}
x^{\prime }(s)\cos{\theta}-y^{\prime }(s)\sin{\theta}=\cos{(2\lambda
s+\theta)}, \qquad
x^{\prime }(s)\sin{\theta}+y^{\prime }(s)\cos{\theta}=\sin{(2\lambda
s+\theta)}, \\
x(s)\cos{\theta}-y(s)\sin{\theta}=-\frac{1}{2\lambda}(\sin{\theta}-\sin{%
(2\lambda s+\theta)}), \\
x(s)\sin{\theta}+y(s)\cos{\theta}=\frac{1}{2\lambda}(\cos{\theta}-\cos{%
(2\lambda s+\theta)}), \qquad
x^2(s)+y^2(s)=\frac{1}{2\lambda^2}(1-\cos{2\lambda s}),
\end{gather*}
we have
\begin{align*} 
Y_{s}={}&\bigl(x^{\prime }(s)\cos \theta-y^{\prime }(s)\sin \theta \bigr)
\mathring e_1+\bigl(x^{\prime }(s)\sin \theta+y^{\prime }(s)\cos \theta \bigr)
\mathring e_2|_{Y(s,\theta)} \\
={}&\cos{(2\lambda s+\theta)}\mathring e_1+\sin{(2\lambda s+\theta)}\mathring
e_2, \\
Y_{\theta}={}&\bigl(x^{\prime }_{1}(\theta)-x(s)\sin \theta-y(s)\cos \theta
\bigr)\mathring e_1+\bigl(x^{\prime }_{2}(\theta)+x(s)\cos \theta-y(s)\sin
\theta \bigr)\mathring e_2 +\bigl(\Theta(\mathcal{C}^{\prime }(\theta))\\
&+2x^{\prime }_{2}(\theta)
\bigl(x(s)\cos \theta-y(s)\sin \theta \bigr)-2x^{\prime }_{1}(\theta)\bigl(
x(s)\sin \theta+y(s)\cos{\theta}\bigr)
+x^2(s)+y^2(s)\bigr)T \\
={}&\left(x^{\prime }_{1}(\theta)-\frac{1}{2\lambda}(\cos{\theta}-\cos{
(2\lambda s+\theta)})\right)\mathring e_1+\left(x^{\prime }_{2}(\theta)-\frac{1
}{2\lambda}(\sin{\theta}-\sin{(2\lambda s+\theta)})\right)\mathring e_2 \\
& +\left(\Theta(\mathcal{C}^{\prime }(\theta))-x^{\prime }_{2}(\theta)
\frac{1}{\lambda}(\sin{\theta}-\sin{(2\lambda s+\theta)})-x^{\prime
}_{1}(\theta)\frac{1}{\lambda}(\cos{\theta}-\cos{(2\lambda s+\theta)}) \right.\\
&\left. +\frac{1}{2\lambda^2}(1-\cos{2\lambda s})\right)T.
\end{align*}
Therefore,
\begin{align}
Y_{s}\wedge Y_{\theta}={}&\left[x^{\prime }_{2}(\theta)\cos{(2\lambda
s+\theta)}-x^{\prime }_{1}(\theta)\sin{(2\lambda s+\theta)}+\frac{\sin{
2\lambda s}}{2\lambda}\right]\mathring e_{1}\wedge \mathring e_{2} \nonumber\\
& +[\cos{(2\lambda s+\theta)}\big<Y_{\theta},T\big>]\mathring
e_{1}\wedge T+[\sin{(2\lambda s+\theta)}\big<Y_{\theta},T\big>]\mathring
e_{2}\wedge T \nonumber\\
={}&\left[A(\theta)\cos{2\lambda s}+\left(\frac{1}{2\lambda}
-B(\theta)\right)\sin{2\lambda s}\right]\mathring e_{1}\wedge \mathring e_{2}\nonumber
\\
&+[\cos{(2\lambda s+\theta)}\big<Y_{\theta},T\big>]\mathring
e_{1}\wedge T+[\sin{(2\lambda s+\theta)}\big<Y_{\theta},T\big>]\mathring
e_{2}\wedge T,\label{cofra1}
\end{align}
where
\begin{gather}
A(\theta)=x^{\prime }_{2}(\theta)\cos{\theta}-x^{\prime }_{1}(\theta)\sin{
\theta},\qquad
B(\theta)=x^{\prime }_{2}(\theta)\sin{\theta}+x^{\prime }_{1}(\theta)\cos{
\theta}, \nonumber\\
\big<Y_{\theta},T\big>=\frac{1}{\lambda}\left[\left(B(\theta)-\frac{1}{
2\lambda}\right)\cos{2\lambda s}+A(\theta)\sin{2\lambda s}+D(\theta)\right],\nonumber
\\
D(\theta)=\lambda \Theta(\mathcal{C}^{\prime }(\theta))+\left(\frac{1}{
2\lambda}-B(\theta)\right). \label{cofra2}
\end{gather}
From \eqref{cofra1}, we conclude that $Y$ is \emph{an immersion} if and
only if either
\[
\left[A(\theta)\cos{2\lambda s}+\left(\frac{1}{2\lambda}-B(\theta)\right)\sin
{2\lambda s}\right]\neq 0\qquad \text{or}\qquad \big<Y_{\theta},T\big>\neq 0.
\]
For the constructed surface $Y$ in \eqref{pacms}, we always assume it is
defined on a region such that $Y$ is an immersion and $\Sigma_{\mathcal{C}}$
is the constant $p$-mean curvature surface defined by such an immersion~$Y$.
A point $p\in \Sigma_{\mathcal{C}}$ is \emph{a singular point} if and only
if $\big<Y_{\theta},T\big>=0$. Thus at a singular point, we must have
\[
\left[A(\theta)\cos{2\lambda s}+\left(\frac{1}{2\lambda}-B(\theta)\right)\sin
{2\lambda s}\right]\neq 0.
\]

Now, we proceed to compute the invariants for $Y$. From the construction of $
Y $, we see that~$(s,\theta)$ is a compatible coordinate system and we are able to
choose the characteristic direction~${e_{1}=Y_{s}}$, and hence
\[
e_{2}=Je_{1}=-\sin{(2\lambda s+\theta)}\mathring e_1+\cos{(2\lambda s+\theta)
}\mathring e_2.
\]

The $\alpha $-function is a function defined on the regular part that
satisfies \[\alpha e_{2}+T=a\sqrt{1+\alpha ^{2}}Y_{s}+b\sqrt{1+\alpha ^{2}}
Y_{\theta }=a\sqrt{1+\alpha ^{2}}e_{1}+b\sqrt{1+\alpha ^{2}}Y_{\theta }\] for
some functions $a$ and $b$. This is equivalent to, comparing the alike terms,
\begin{gather*}
-\alpha \sin {(2\lambda s+\theta )} =a\sqrt{1+\alpha ^{2}}\cos {(2\lambda
s+\theta )}\\
\phantom{-\alpha \sin {(2\lambda s+\theta )} =}{}+b\sqrt{1+\alpha ^{2}}\left( x_{1}^{\prime }(\theta )-\frac{1}{
2\lambda }(\cos {\theta }-\cos {(2\lambda s+\theta )})\right), \\
\alpha \cos {(2\lambda s+\theta )}=a\sqrt{1+\alpha ^{2}}\sin {(2\lambda
s+\theta )}\\
\phantom{\alpha \cos {(2\lambda s+\theta )}=}{}+b\sqrt{1+\alpha ^{2}}\left( x_{2}^{\prime }(\theta )-\frac{1}{
2\lambda }(\sin {\theta }-\sin {(2\lambda s+\theta )})\right), \\
1 =b\sqrt{1+\alpha ^{2}}\big<Y_{\theta },T\big>.
\end{gather*}

We thus have
\begin{gather}
a =\frac{-2\lambda \bigl( x^{\prime }_{1}( \theta ) \cos (
2\lambda s+\theta ) +x^{\prime }_{2}( \theta ) \sin (
2\lambda s+\theta ) \bigr) -( 1-\cos 2\lambda s) }{
2\lambda \sqrt{1+\alpha ^{2}} \langle Y_{\theta },T \rangle },\nonumber \\
b =\frac{1}{\sqrt{1+\alpha ^{2}}\big<Y_{\theta },T\big>}, \nonumber\\
\alpha =\frac{x_{2}^{\prime }(\theta )\cos {(2\lambda s+\theta )}
-x_{1}^{\prime }(\theta )\sin {(2\lambda s+\theta )}+\frac{\sin {2\lambda s}
}{2\lambda }}{\big<Y_{\theta },T\big>} \nonumber\\
\phantom{\alpha}{} =\lambda \frac{A(\theta )\cos {2\lambda s}+\bigl( \frac{1}{2\lambda }
-B(\theta )\bigr) \sin {2\lambda s}}{\bigl( B(\theta )-\frac{1}{2\lambda }
\bigr) \cos {2\lambda s}+A(\theta )\sin {2\lambda s}+D(\theta )}.\label{forinv}
\end{gather}
Let $V=\bigl( A(\theta ),\frac{1}{2\lambda }-B(\theta )\bigr) $ and \smash{$
\Vert V \Vert =\sqrt{[ A( \theta ) ] ^{2}+
\bigl[ \frac{1}{2\lambda }-B(\theta )\bigr] ^{2}}$}.
If $V=0$, then $\alpha=0$. If $V\neq 0$, then we can write
\smash{$
\frac{V}{ \Vert V \Vert }=(\sin {\zeta (\theta )},\cos {\zeta
(\theta )})$},
for some function $\zeta {(\theta )}$. The functions $\alpha $, $a$, and $b$
can be further written as
\begin{gather*}
\alpha =\lambda \frac{\sin {\zeta (\theta )}\cos {2\lambda s}+\cos {\zeta
(\theta )}\sin {2\lambda s}}{-\cos {\zeta (\theta )}\cos {2\lambda s}+\sin {
\zeta (\theta )}\sin {2\lambda s}+\frac{D(\theta )}{ \Vert V
\Vert }} =\lambda \frac{\sin {\bigl(2\lambda s+\zeta {(\theta )}\bigr)}}{G(\theta
)-\cos {\bigl(2\lambda s+\zeta {(\theta )}\bigr)}}, \\
a=\frac{-2\lambda [ \sin \zeta ( \theta ) \sin (
2\lambda s) -\cos \zeta ( \theta ) \cos 2\lambda s] -
\frac{1}{\Vert V \Vert }}{2\sqrt{1+\alpha ^{2}}\bigl[ \frac{
D( \theta ) }{\Vert V \Vert }-\cos ( 2\lambda
s+\zeta ( \theta ) ) \bigr] } \\
\phantom{a}{} =\frac{2\lambda \Vert V \Vert \cos ( 2\lambda s+\zeta
( \theta ) ) -1}{2\Vert V \Vert \sqrt{1+\alpha
^{2}}[ G( \theta ) -\cos ( 2\lambda s+\zeta (
\theta ) ) ] }, \\
b =\frac{\frac{\lambda }{\Vert V \Vert }}{\sqrt{1+\alpha ^{2}}
\bigl[ \frac{D( \theta ) }{\Vert V \Vert }-\cos
( 2\lambda s+\zeta ( \theta ) ) \bigr] } =\frac{\frac{\lambda }{\Vert V \Vert }}{\sqrt{1+\alpha ^{2}}
[ G( \theta ) -\cos ( 2\lambda s+\zeta ( \theta
) ) ] },
\end{gather*}
where
\begin{equation}\label{gdefn}
G(\theta )=\frac{D(\theta )}{\Vert V \Vert }=\frac{D(\theta )}{
\sqrt{(A(\theta ))^{2}+\bigl(\frac{1}{2\lambda }-B(\theta )\bigr)^{2}}}.
\end{equation}

Next, we normalize the three invariants $\alpha $, $a$, and $b$.
Firstly, we choose another compatible coordinates ($ \tilde{s}=s+\Gamma
( \theta )$, $\tilde{\theta}=\Psi ( \theta ) $), for some $\Gamma ( \theta ) $ and $\Psi ( \theta
) $. From the transformation law of the induced metric
$
\tilde{a}=a+b\Gamma^{\prime }( \theta )$, $\tilde{b}=b\Psi
^{\prime }( \theta )$,
this can be chosen so that
\begin{gather*}
\Gamma ^{\prime}( \theta ) =\frac{-2\lambda \Vert V(
\theta ) \Vert G( \theta ) +1}{2\lambda } =\frac{1}{2\lambda }-D( \theta ) , \notag
\end{gather*}
or equivalently,
\begin{equation*}
\Gamma ( \theta ) =\frac{\theta }{2\lambda }-\int D( \theta
) {\rm d}\theta .
\end{equation*}
If we further choose $\Psi $ such that $\tilde{\theta}=\Psi ( \theta
) =2\lambda \int \Vert V( \theta ) \Vert
{\rm d}\theta $, then in terms of the compatible coordinates $\bigl( \tilde{s},\tilde{\theta}\bigr) $, the three invariants read
\begin{gather*}
\tilde{a} = \frac{-\lambda }{\sqrt{1+\tilde{\alpha}^{2}}}, \qquad
\tilde{b} = \frac{2\lambda ^{2}}{\sqrt{1+\tilde{\alpha}^{2}}\big\{
G\bigl( \Psi ^{-1}\bigl( \tilde{\theta}\bigr) \bigr) -\cos \bigl[ 2\lambda
\tilde{s}-2\lambda \Gamma \bigl( \Psi ^{-1}\bigl( \tilde{\theta}\bigr)
\bigr) +\zeta \bigl( \Psi ^{-1}\bigl( \tilde{\theta}\bigr) \bigr) \bigr]
\big\} }, \\
\tilde{\alpha} = \lambda \frac{\sin \bigl( 2\lambda \tilde{s}-2\lambda
\Gamma ( \theta ) +\zeta \bigl( \Psi ^{-1}\bigl( \tilde{\theta}
\bigr) \bigr) \bigr) }{G\bigl( \Psi ^{-1}\bigl( \tilde{\theta}\bigr)
\bigr) -\cos \bigl[ 2\lambda \tilde{s}-2\lambda \Gamma \bigl( \Psi
^{-1}\bigl( \tilde{\theta}\bigr) \bigr) +\zeta \bigl( \Psi ^{-1}\bigl(
\tilde{\theta}\bigr) \bigr) \bigr] },
\end{gather*}
where $D$ and $G$ are defined in \eqref{cofra2} and \eqref{gdefn}, respectively.

We summarize the above discussion as a theorem in the following.

\begin{Theorem}
\label{construction} The coordinate system $(s,\theta)$ for $Y$ in \eqref{pacms}
is compatible. If $V=0$, then $\alpha=0$. If $V\neq 0$, then the new coordinate system $\bigl(\tilde{s},\tilde{\theta}\bigr)$, where
$
\tilde{s}=s+\Gamma(\theta)$, $ \tilde{\theta}=\Psi( \theta)$,
with
\[
\Gamma(\theta)=\frac{\theta }{2\lambda}-\int D( \theta) {\rm d}\theta,\qquad
\Psi(\theta)=2\lambda \int \Vert V(\theta) \Vert {\rm d}\theta ,
\]
is normal. In terms of the
normal coordinates, the invariants of $Y$ are given by
\begin{equation}
\zeta _{1}\bigl( \tilde{\theta}\bigr) = \zeta \bigl( \Psi ^{-1}\bigl(
\tilde{\theta}\bigr) \bigr) -2\lambda \Gamma \bigl( \Psi ^{-1}\bigl(
\tilde{\theta}\bigr) \bigr) , \qquad
\zeta _{2}\bigl( \tilde{\theta}\bigr) = G\bigl( \Psi ^{-1}\bigl( \tilde{
\theta}\bigr) \bigr). \label{invariants}
\end{equation}
\end{Theorem}

Particularly, in order to have constant $\zeta_1\bigl(\tilde \theta\bigr)$ and nonzero
constant $\zeta_2\bigl(\tilde \theta\bigr)$, Theorem \ref{construction} suggests the
constant $p$-mean curvature surfaces deformed by curves
\begin{equation} \label{deformedcurve}
\mathcal{C}(\theta)=(x_{1}(\theta),x_{2}(\theta),x_{3}(\theta))=\left(\frac{r}{
\lambda}\sin{\theta},-\frac{r}{\lambda}\cos{\theta},\frac{r(1-r)}{\lambda^2}
\theta\right),
\end{equation}
where $r\neq \frac{1}{2}$. More precisely, we have the following proposition.

\begin{Proposition}
\label{c2value} For any curve $\mathcal{C}(\theta)$ defined as
\eqref{deformedcurve}, the deformed surface $Y(s,\theta)$ has both constant
invariants $\zeta_1\bigl(\tilde \theta\bigr)$ and $\zeta_2\bigl(\tilde \theta\bigr)\neq 0$.
\end{Proposition}

\begin{proof}
We argue by assuming $\zeta_2\bigl(\tilde \theta\bigr)=\zeta_2$ is a constant, $
x_{1}(\theta)=\frac{r}{\lambda}\sin{\theta}$ , and $x_{2}(\theta)=-\frac{r}{
\lambda}\cos{\theta}$ for any $r\neq \frac{1}{2}$. Then \eqref{cofra2}
implies $A(\theta)=0$, $B(\theta)=\frac{r}{\lambda}$, which leads to
\smash{$
\Vert V \Vert =\frac{|1-2r|}{2\lambda}$}.
The second equation of \eqref{invariants} shows that $D(\theta)=\zeta_2
\Vert V \Vert $, and hence
\begin{align*}
\zeta_1\bigl(\tilde \theta\bigr)&= \sin^{-1}\left(\frac{A(\theta)}{\Vert V
\Vert }\right)-2\lambda \left(\frac{\theta}{2\lambda}-\int D(\theta) {\rm d}\theta
\right) =-\theta+\int \zeta_2|1-2r|{\rm d}\theta \\
&=(\zeta_2|1-2r|-1)\theta+\text{\rm const}.
\end{align*}
In order to have $\zeta_1\bigl(\tilde \theta\bigr)$ being constant, we must have
$
\zeta_2|1-2r|=1$.
It is clear to see that $\zeta_2\neq 0$ and $r\neq 0$. The system
\eqref{invariants} immediately shows $D(\theta)=\frac{1}{2\lambda}$, which
gives
\smash{$
x_3^{\prime }(\theta)=\frac{r(1-r)}{\lambda^2}$},
by \eqref{cofra2}. Namely, \smash{$x_3(\theta)=\frac{r(1-r)}{\lambda^2}
\theta+\text{\rm const}$}.

Moreover, the new coordinates can be obtained by
\[
\tilde \theta=\Psi(\theta)=|1-2r|\theta+\text{\rm const}\qquad \text{and} \qquad
\Gamma(\theta),\]
up to a constant.
\end{proof}

\section{Examples}\label{examples}

It is easy to see that the Pansu sphere can be obtained by deforming the following curves
\[
\mathcal{C}_1(\theta)=(0,0,\text{\rm const}) \qquad\text{ or }\qquad \mathcal{C}_2(\theta)=\left(\frac{1
}{\lambda}\sin{\theta},-\frac{1}{\lambda}\cos{\theta},\text{\rm const}\right).
\]

Using a similar idea as Theorem \ref{construction} and Proposition \ref{c2value}, we obtain curves $\mathcal{C}(\theta)$ that result in constant $p$-mean curvature surfaces with
constant $\zeta_2$ and linear $\zeta_1\bigl(\tilde \theta\bigr)$ in Sections~\ref{section5.1} and~\ref{section5.3}. We collect $\mathcal C(\theta)$ in Tables~\ref{tab1} and \ref{tab2} as follows.

\begin{table}[ht]\centering\renewcommand{\arraystretch}{1.2}\vspace{-2mm}
\caption{Examples of $\mathcal C(\theta)$ for constant $p$-mean curvature surfaces.}\label{tab1}\vspace{1mm}

\begin{tabular}{|c|c|c|}
\hline
$\mathcal C(\theta)$ & constant $\zeta_1$ & linear $\zeta_1$ \\ \hline
$\zeta_2>1$ & \begin{tabular}[c]{@{}c@{}}$\bigl(\frac{r}{\lambda}\sin\theta,-\frac{r}{\lambda}\cos\theta,\frac{r(1-r)}{\lambda^2}\theta\bigr)$\\ $0<r<1, r\neq \frac{1}{2}$ \\ $\zeta_2=\frac{1}{|1-2r|}$\end{tabular} & \multirow{4}{*}{\begin{tabular}[c]{@{}c@{}}$m\neq 1\colon (x_1(\theta),x_2(\theta),x_3(\theta))$\\ $x_1(\theta)=\frac{\sin\theta}{2\lambda}-\frac{\sin((m-1)\theta)}{2\lambda k(m-1)}$\\ $x_2(\theta)=-\frac{\cos\theta}{2\lambda}-\frac{\cos((m-1)\theta)}{2\lambda k(m-1)}$\\ $x_3(\theta)=\frac{1+k^2(m-1)}{4\lambda^2k^2(m-1)}\theta-\frac{\sin(m\theta)}{4\lambda^2 k(m-1)}$ \\ $m=1\colon \bigl(\frac{\sin\theta}{2\lambda}-\frac{\theta}{2\lambda k},-\frac{\cos\theta}{2\lambda},\frac{k\theta-\theta\cos\theta}{4\lambda^2k}\bigr)$\\ $\zeta_1=m\theta+{\rm const}$ and $\zeta_2=k>0$\end{tabular}} \\ \cline{1-2}
\multirow{2}{*}{$\zeta_2=1$} & \multirow{2}{*}{Pansu sphere} & \\
 & & \\ \cline{1-2}
$0<\zeta_2<1$ & \begin{tabular}[c]{@{}c@{}}$\bigl(\frac{r}{\lambda}\sin\theta,-\frac{r}{\lambda}\cos\theta,\frac{r(1-r)}{\lambda^2}\theta\bigr)$\\ $r<0$ or $r>1$ \\ $\zeta_2=\frac{1}{|1-2r|}$\end{tabular} & \\ \hline
$\zeta_2=0$ & \begin{tabular}[c]{@{}c@{}}$\bigl(\frac{\beta}{4\lambda},0,\frac{\beta}{4\lambda}-\frac{\theta}{4\lambda^2}\bigr)$\\ $\beta=\ln |\sec\theta+\tan\theta|$\\ $\zeta_1=0$\end{tabular} & \begin{tabular}[c]{@{}c@{}}$\bigl(\cos\theta,\sin\theta,-\bigl(\theta+\frac{\theta}{2\lambda^2}\bigr)\bigr)$ \\ $\zeta_1=-\theta+{\rm const}$\end{tabular} \\ \hline
\end{tabular}\vspace{-4mm}
\end{table}

\begin{table}[ht]\centering\renewcommand{\arraystretch}{1.2}
\caption{Examples of $\mathcal C(\theta)$ for $p$-minimal surfaces.}\label{tab2}\vspace{1mm}

\begin{tabular}{|c|cc|}
\hline
$\mathcal C(\theta)$ & \multicolumn{1}{c|}{constant $\zeta_1$} & linear $\zeta_1$ \\ \hline
$\zeta_2>0$ & \multicolumn{1}{c|}{type I} & \begin{tabular}[c]{@{}c@{}}$(r\sin\theta,-r\cos\theta, z(\theta))$\\ $z'(\theta)+r^2>0$\\ $\zeta_1=-r\theta$\end{tabular} \\ \hline
$\zeta_2<0$ & \multicolumn{1}{c|}{type II, III} & \begin{tabular}[c]{@{}c@{}}$(r\sin\theta,-r\cos\theta, z(\theta))$\\ $z'(\theta)+r^2<0$\\ $\zeta_1=-r\theta$\end{tabular} \\ \hline
special type I & \multicolumn{2}{c|}{\begin{tabular}[c]{@{}c@{}}degenerate case:$\bigl(-\theta,0,\frac{\sin(2\theta)-2\theta}{4}\bigr)$\\ or\\ entire graph: $u=0$\end{tabular}} \\ \hline
special type II & \multicolumn{2}{c|}{$u=xy+g(y)$} \\ \hline
\end{tabular}\vspace{-3mm}
\end{table}

\subsection[Examples of constant p-mean curvature surfaces]{Examples of constant $\boldsymbol{p}$-mean curvature surfaces}\label{section5.1}

\begin{Proposition}\label{c1linear}\samepage Given any curve
\[
\mathcal{C}(\theta)=\left(\frac{1}{\lambda}\sin{\theta},-\frac{1}{\lambda}\cos{
\theta},\frac{k-1}{2\lambda^2}\theta+\text{\rm const}\right),
\]
the deformed surface $Y(s,\theta)$ has the invariants $\zeta_1\bigl(\tilde
\theta\bigr)=(k-1)\tilde \theta +\text{\rm const}$ and $\zeta_2\bigl(\tilde \theta\bigr)=k$, where~${k
\in \R}$.
\end{Proposition}

\begin{Remark}
It is easy to see that the surfaces obtained by curves given in Proposition \ref{c1linear} are not rotationally symmetric since $\zeta_1^{\prime
}=k-1<\zeta_2$ by \eqref{coinro1}.
\end{Remark}

\begin{Proposition}\label{indepc1}
For any constant $k>0$ and $m$, there exist constant $p$-mean curvature surfaces~$Y(s,\theta)$ defined as \eqref{pacms} with invariants $\zeta_1(\theta)=m\theta+\text{\rm const}$ and $\zeta_2=k$.
\end{Proposition}
\begin{proof}
It suffices to solve the system \eqref{invariants}. In order to obtain a surface with linear $\zeta_1(\theta)=m\theta$ for any given nonzero constant $\zeta_2=k$, we assume
\begin{equation}\label{ae}
A(\theta)=\frac{1}{2\lambda k}\sin(m\theta)\qquad\text{and}\qquad \frac{1}{2\lambda}-B(\theta)=\frac{1}{2\lambda k}\cos(m\theta).
\end{equation}
It results in $\|V\|=\frac{1}{2\lambda k}$, and
\begin{align*}
\zeta_1(\theta)&=\sin^{-1}\left(\frac{A(\theta)}{\|V\|}\right)-2\lambda\left(\frac{\theta}{2\lambda}-\int D(\theta)~{\rm d}\theta\right)=m\theta-\theta+2\lambda\int \zeta_2(\theta)\|V\|~{\rm d}\theta\\
&=m\theta-\theta+2\lambda k\int \frac{1}{2\lambda k}~{\rm d}\theta=m\theta+\text{\rm const}.
\end{align*}
Next we solve for $x_1'(\theta)$ and $x_2'(\theta)$ from \eqref{cofra2}, that is,
\[
\frac{1}{2\lambda k}\sin(m\theta)=x^{\prime }_{2}(\theta)\cos{\theta}-x^{\prime }_{1}(\theta)\sin{
\theta}, \qquad
\frac{1}{2\lambda}-\frac{1}{2\lambda k}\cos(m\theta)=x^{\prime }_{2}(\theta)\sin{\theta}+x^{\prime }_{1}(\theta)\cos{
\theta}.
\]
It is easy to see that
\[
x_1'(\theta)=\frac{1}{2\lambda}\cos(\theta)-\frac{1}{2\lambda k}\cos((m-1)\theta), \qquad x_2'(\theta)=\frac{1}{2\lambda}\sin(\theta)+\frac{1}{2\lambda k}\sin((m-1)\theta),\]
and hence for $m\neq 1$,
\begin{align}
x_1(\theta)&=\frac{1}{2\lambda}\sin(\theta)-\frac{1}{2\lambda k(m-1)}\sin((m-1)\theta)+\text{\rm const},\nonumber\\
x_2(\theta)&=-\frac{1}{2\lambda}\cos(\theta)-\frac{1}{2\lambda k(m-1)}\cos((m-1)\theta)+\text{\rm const}.\label{ad1}
\end{align}
The equation \eqref{cofra2} also suggests
\[x_3'(\theta)=\frac{1+k^2(m-1)}{4\lambda^2k^2(m-1)}-\frac{m\cos(m\theta)}{4\lambda^2k(m-1)},\]
and then we have
\begin{equation}\label{ad2}
x_3(\theta)=\frac{1+k^2(m-1)}{4\lambda^2k^2(m-1)}\theta-\frac{\sin(m\theta)}{4\lambda^2k(m-1)}+\text{\rm const}.
\end{equation}
Therefore, deforming such curves $\mathcal{C}(\theta)=(x_{1}(\theta),x_{2}(\theta),x_{3}(\theta))$ defined by \eqref{ad1} and \eqref{ad2} gives surfaces with nonzero $\zeta_2=k$ and linear $\zeta_1(\theta)=m\theta+\text{\rm const}$ for all $m \neq 1$.

When $m =1$, direct computations from \eqref{ae} imply
\begin{gather*}
x_1(\theta)=\frac{1}{2\lambda}\sin(\theta)-\frac{\theta}{2\lambda k}+\text{\rm const},\qquad
x_2(\theta)=-\frac{1}{2\lambda}\cos(\theta)+\text{\rm const},\\
x_3(\theta)=\frac{k\theta-\theta\cos\theta}{4\lambda^2k}.\tag*{\qed}
\end{gather*} \renewcommand{\qed}{}
\end{proof}

\begin{Example}
If
\begin{align*}
C(\theta )& =(x_{1}(\theta ),x_{2}(\theta ),x_{3}(\theta )) \\
& =\left(\frac{1}{4\lambda }\ln {|\sec {\theta }+\tan {\theta }|}+c_{3},c_{4},
\frac{c_{5}}{4\lambda }\ln {|\sec {\theta }+\tan {\theta }|}-\frac{1}{
4\lambda ^{2}}\theta +c_{6}\right),
\end{align*}
then $\zeta _{1}\bigl(\tilde{\theta}\bigr)=0$, $\zeta _{2}\bigl(\tilde{\theta}\bigr)=0.$
\end{Example}

\subsection{Basic properties of surfaces of special type I}\label{section5.2}
For $p$-minimal surfaces of special type I, we have the first fundamental form, in terms of normal coordinates $(x,y)$,
\[I={\rm d}x\otimes {\rm d}x+\left(\frac{1+\alpha^{2}}{\alpha^4}\right){\rm d}y\otimes {\rm d}y,\]
so that $I$ degenerates along the curve where $\alpha$ blows up. Recall that the parametrization of the surface $Y$ is
\[Y(r,\theta)=(x(\theta)+r\cos{\theta},y(\theta)+r\sin{\theta},z(\theta)+ry(\theta)\cos{\theta}-rx(\theta)\sin{\theta}).\]
We have
\[
Y_{r}=(\cos{\theta},\sin{\theta},y(\theta)\cos{\theta}-x(\theta)\sin{\theta}),\qquad
Y_{\theta}=\bigl(x'(\theta)-r\sin{\theta},y'(\theta)+r\cos{\theta},*\bigr),
\]
where
\[*=z'(\theta)+ry'(\theta)\cos{\theta}-y(\theta)\sin{\theta}-x'(\theta)\sin{\theta}-x(\theta)\cos{\theta}.\]
Then
\begin{align*}
Y_{r}\times Y_{\theta}={}&\left|\begin{matrix}
i&j&k\\
\cos{\theta}&\sin{\theta}&y(\theta)\cos{\theta}-x(\theta)\sin{\theta}\\
x'(\theta)-r\sin{\theta}&y'(\theta)+r\cos{\theta}&*
\end{matrix}\right|\\
={}&\rho\bigl(\sin{\theta}\bigl(y'(\theta)\cos{\theta}-x'(\theta)\sin{\theta}\bigr)-y(\theta),\\
&-\cos{\theta}\bigl(y'(\theta)\cos{\theta}-x'(\theta)\sin{\theta}\bigr)+x(\theta),1\bigr),
\end{align*}
where
$\rho=r+\bigl(y'(\theta)\cos{\theta}-x'(\theta)\sin{\theta}\bigr)$.
For $p$-minimal surfaces of special type I, we have
 \[
 \alpha=\frac{1}{r+\bigl(y'(\theta)\cos{\theta}-x'(\theta)\sin{\theta}\bigr)}.
 \]
Therefore, $Y_{r}$ and $Y_{\theta}$ are linearly dependent along the curve when $\alpha$ blows up.

For constant $p$-mean curvature surfaces of special type I, \eqref{cofra1} and \eqref{cofra2} immediately imply that $Y_{s}\wedge Y_{\theta}=0$ if and only if
\begin{gather*}
0=A(\theta)\cos{2\lambda s}+\left(\frac{1}{2\lambda}-B(\theta)\right)\sin{2\lambda s},\\
-D(\theta)=-\left(\frac{1}{2\lambda}-B(\theta)\right)\cos{2\lambda s}+A(\theta)\sin{2\lambda s},
\end{gather*}
that is,
\[\left(\begin{matrix}
\cos{2\lambda s}&\sin{2\lambda s}\\\sin{2\lambda s}&-\cos{2\lambda s}
\end{matrix}\right)\left(\begin{matrix}A(\theta)/\|V\| \\ \bigl(\frac{1}{2\lambda}-B(\theta)\bigr)/\|V\|\end{matrix}\right)=\left(\begin{matrix}
0\\-G(\theta)\end{matrix}\right).\]
This implies that $Y_{s}\wedge Y_{\theta}=0$ holds if and only if $G(\theta)=\pm 1$, namely, it happens only on the surface $Y$ of special type I at points where the function $\alpha$ blows up.

\subsection[Examples of p-minimal surfaces]{Examples of $\boldsymbol{p}$-minimal surfaces}\label{section5.3}

In what follows, we give some $p$-minimal surfaces of special type~II (i.e., $\zeta_2<0$ and linear $\zeta_1$). We first recall in \cite{ChiuH/LiuH:2022} that $\mathcal C(\theta)=(x(\theta),y(\theta),z(\theta))$ satisfying
\begin{align*}
\zeta_1(\theta)&=-\Gamma(\theta)+y'(\theta)\cos\theta-x'\sin\theta,\\
\zeta_2(\theta)&=z'(\theta)+x(\theta)y'(\theta)-y(\theta)x'(\theta)-\bigl(y'(\theta)\cos\theta-x'(\theta)\sin\theta\bigr)^2,
\end{align*}
where $\Gamma(\theta)=\int x'(\theta)\cos\theta+y'(\theta)\sin\theta {\rm d}\theta$, will result in a $p$-minimal surface. For any nonzero $r\in\R$, if $x(\theta)=r\sin\theta$ and $y(\theta)=-r\cos\theta$,
then $\Gamma(\theta)=r\theta$ (up to a constant), $\zeta_1(\theta)=-r\theta$, and~${\zeta_2(\theta)=z'(\theta)+r^2}$.
We choose $z(\theta)$ such that $z'(\theta)+r^2<0$ to have negative $\zeta_2$.

\subsection*{Acknowledgements}

The first author's research was supported in part by NSTC 112-2115-M-007-009-MY3. The second author's research was supported in part by NSTC 110-2115-M-167-002-MY2 and NSTC 112-2115-M-167-002-MY2. The third author's research was supported in part by NSTC 112-2628-M-032-001-MY4. We all thank the anonymous referees for carefully reading our manuscript and their insightful comments and suggestions for improving the article.

\pdfbookmark[1]{References}{ref}
\LastPageEnding

\end{document}